\magnification=\magstep1

\def\refs{\medskip\hangindent=25pt\hangafter=1\noindent}

\font\sc=cmcsc10

\centerline {ON THE BEHAVIOR OF THE COVARIANCE MATRICES IN A}
\centerline {MULTIVARIATE CENTRAL LIMIT THEOREM}
\centerline {UNDER SOME MIXING CONDITIONS}
\bigskip
\bigskip

\noindent Richard C.\ Bradley \hfil\break
Department of Mathematics \hfil\break
Indiana University \hfil\break
Bloomington \hfil\break
Indiana 47405 \hfil\break
USA \hfil\break

\noindent bradleyr@indiana.edu \hfil\break
\vskip 0.5 in

   {\bf Abstract.} In a paper that appeared in 2010, 
C.\ Tone proved a multivariate central limit theorem for 
some strictly stationary random fields of random vectors satisfying certain mixing conditions.
The ``normalization'' of a given ``partial sum'' 
(or ``block sum'') involved matrix multiplication by 
a ``standard $-1/2$ power'' of its covariance matrix (a
symmetric, positive definite matrix), and the limiting
multivariate normal distribution had the identity matrix
as its covariance matrix.  
The mixing assumptions in Tone's result implicitly
imposed an upper bound on the ratios of the largest to
the smallest eigenvalues in the covariance matrices of the
partial sums.
The purpose of this note is to show that in Tone's result,
for the entire collection of the covariance matrices of
the partial sums, there is essentially no other 
restriction on the relative magnitudes of the eigenvalues
or on the (orthogonal) directions of the 
corresponding eigenvectors.
For simplicity, the example given in this note will
involve just random sequences, not the broader context
of random fields.
\hfil\break

\hfil\break

\hfil\break
 
\vfill\eject

   {\bf 1. Introduction.}  A multivariate central limit theorem was proved by C.\ Tone [26] for some strictly
stationary random fields of random vectors satisfying 
certain mixing conditions.
As in a somewhat related result in [6] under 
different dependence assumptions, the ``normalization'' 
of a given ``partial sum'' (or ``block sum'')
involved matrix multiplication by 
a ``standard $-1/2$ power'' of its covariance matrix (a
symmetric, positive definite matrix), and the limiting
multivariate normal distribution had the identity matrix
as its covariance matrix.  (More on that below.)\ \ 
The mixing assumptions in Tone's [26] result implicitly
imposed an upper bound on the ratios of the largest to
the smallest eigenvalues in the covariance matrices of the
partial sums.
The purpose of this note is to show that in Tone's result,
for the entire collection of the covariance matrices of
the partial sums, there is essentially no other 
restriction on the relative magnitudes of the eigenvalues
or on the (orthogonal) directions of the 
corresponding eigenvectors.
This will be elucidated with an example described in
Theorem 1.4 below, after a special case of Tone's result
is stated in Theorem 1.3.
For simplicity, our attention in this note will be confined
to just {\it sequences\/} (of random vectors), instead of 
the broader context of random fields. 

   First, Sections 1.1 and 1.2 will give some definitions
and notations and will also briefly review some well
known, standard, elementary mathematics that will be 
needed.
\hfil\break

   {\sc Notations 1.1.} In what follows, the entries
of matrices are real numbers. 
The transpose of any given matrix $M$ will be 
denoted $M^t$.

   Now suppose $m$ is a positive integer.
In some of the notations below, the dependence on this
given positive integer $m$ will be tacitly understood
and not indicated explicitly. \smallskip
   
   (A) A given element $x \in {\bf R}^m$ will be represented
as a ``column vector'' (an $m \times 1$ matrix):
$x := [x_1, x_2, \dots, x_m]^t$.
For such an $x$, denote the Euclidean norm as
$\|x\| := (x_1^2 + x_2^2 + \dots + x_m^2)^{1/2}$.
The origin in ${\bf R}^m$ will be denoted 
${\bf 0}_m := [0, 0, \dots, 0]^t$.
\smallskip

   (B) A symmetric $m \times m$ matrix $A$ is
``positive semi-definite'' if $x^tAx \geq 0$ for
all $x \in {\bf R}^m$, and $A$ is ``positive definite''
if $x^tAx > 0$ (strict inequality) for all
$x \in {\bf R}^m - \{{\bf 0}_m\}$.
\smallskip 
   
   (C) If $A$ is a symmetric, positive definite
(hence nonsingular)
$m \times m$ matrix and $r$ is a real number, then
$A^r$ denotes the symmetric, positive definite
``$r^{\rm th}$ power'' matrix of $A$.

(It is of course defined by $A^r := UD^rU^t$ where
(i) $U$ is an ($m \times m$) orthogonal matrix and $D$ a
diagonal matrix such that 
$A = U D U^t$ and 
(ii) $D^r$ is the diagonal matrix in which,
for each $i \in \{1, \dots, m\}$, the $i^{th}$ diagonal
element is $d_i^r$ where $d_i$ (a positive number,
an eigenvalue of $A$)
is the $i^{th}$ diagonal element of $D$.
The matrix $A^r$ will thereby be uniquely defined,
even though in general the choice of matrices $U$ and $D$
in this procedure is not unique.)
\smallskip

   (D) For any given symmetric, positive definite 
$m \times m$ matrix $A = (a_{ij}, 1 \leq i,j \leq m)$, 
define the following two quantities:
$$ \eqalignno { 
\eta_{\rm min}(A) &:= 
\min_{x \in {\bf R}\uparrow m: \|x\| = 1} x^tAx, 
\quad {\rm and}
& (1.1) \cr
\eta_{\rm max}(A) &:= 
\max_{x \in {\bf R}\uparrow m: \|x\| = 1} x^tAx.
& (1.2) \cr
} $$
(As in (1.1) and (1.2), a notation of
the form $\alpha^\beta$ in a subscript or superscript will 
typically be written as $\alpha \uparrow \beta$ for typographical convenience.)\ \ 
In (1.1)-(1.2), the min and max are both achieved for 
elements $x$ on the unit sphere, and they are equal respectively to the smallest and largest eigenvalues 
of $A$.
Each entry $a_{ij}$ of $A$ satisfies 
$|a_{ij}| \leq \eta_{\rm max}(A)$. 

\smallskip

   (E) For any two positive numbers $a$ and $b$ such that
$a < b$, let $\Lambda_{(m,a,b)}$ denote the set of all
symmetric, positive definite $m \times m$ matrices $A$
such that 
$a \leq \eta_{\rm min}(A) \leq \eta_{\rm max}(A) \leq b$
(that is, the set of all such matrices whose eigenvalues
are all between $a$ and $b$ inclusive).
\smallskip

   (F) For each $\varepsilon > 0$, let 
${\bf B}_{\rm sym}^{(m)}[\varepsilon]$
denote the set of all symmetric (not necessarily 
positive semi-definite) $m \times m$ matrices 
$B := (b_{ij}, 1 \leq i,j \leq m)$ such that
$|b_{ij}| \leq \varepsilon$ for all 
$(i,j) \in \{1, \dots, m\}^2$.
\smallskip

   (G) If $a$, $b$, and $\varepsilon$ are positive numbers such that $m\varepsilon < a < b$, and 
$A \in \Lambda_{(m,a,b)}$
and $B \in {\bf B}^{(m)}_{\rm sym}[\varepsilon]$, then
$A + B \in \Lambda_{(m,a-m\varepsilon,b+m\varepsilon)}$.    
(The point is that for such a $B$, if $x \in {\bf R}^m$
is such that $\|x\| = 1$, then $|x^tBx| \leq m\varepsilon$
simply by persistent trivial applications of the Cauchy
inequality $|y^tz| \leq \|y\| \cdot \|z\|$ for
$y,z \in {\bf R}^m$.)
\hfil\break

   {\sc Notations 1.2.} Now suppose $(\Omega, {\cal F}, P)$ 
is a probability space.
Again suppose $m$ is a positive integer.
\smallskip

   (A) An ``${\bf R}^m$-valued random variable'' is a 
random vector with $m$ (random real) coordinates.
Such random vectors $V$ will be represented as ``random
column vectors'' (i.e.\ $m \times 1$ random matrices):
$V := [V_1, V_2, \dots, V_m]^t$.
  
   In the case where $E\|V\|^2 < \infty$ (that is, 
$EV_i^2 < \infty$ for each 
$i \in \{1, \dots, m\}$ --- recall Notations 1.1(A)), 
the ($m \times m$) covariance matrix of $V$ will be 
denoted $\Sigma_V$.
If also $EV = {\bf 0}_m$ (that is, $EV_i = 0$ for 
each $i$), then one has the trivial representation 
$\Sigma_V = EVV^t$.
The matrix $\Sigma_V$ is of course (symmetric and)
positive semi-definite. (In the mean ${\bf 0}_m$
case, recall that for any  $x \in {\bf R}^m$,
$x^t\Sigma_Vx = E(x^tV)(x^tV)^t = E(x^tV)^2 \geq 0$).
\smallskip

   (B) Suppose $X := (X_k, k \in {\bf Z})$ is a strictly
stationary sequence of ${\bf R}^m$-valued random variables.
For each $n \in {\bf N}$, define the partial sum
(again, a ``random $m \times 1$ column vector'')
$S_n = S(X,n) := X_1 + X_2 + \dots + X_n$.
(Here and below, ${\bf N}$ denotes the set of all
positive integers.)

   Our work will involve the case where
$EX_0 = {\bf 0}_m$ and $E\|X_0\|^2 < \infty$.
For typographical convenience, 
the covariance matrix of $X_0$ will be written 
$\Sigma_{X(0)}$, and for each $n \in {\bf N}$, the 
covariance matrix of the normalized
partial sum $n^{-1/2}S_n$ will be written 
(with perhaps slight abuse of notation) as
$\Sigma_{S(X,n)/\sqrt n}$ (it is of course equal to
$n^{-1} \Sigma_{S(X,n)}$).
\smallskip

   (C) Next let us turn to measures of dependence.
For any two $\sigma$-fields ${\cal A}$ and ${\cal B}$
($\subset {\cal F}$), define the following four
measures of dependence:  First, define
$$ \alpha({\cal A}, {\cal B}) :=
\sup_{A \in {\cal A}, B \in {\cal B}}
|P(A \cap B) - P(A)P(B)|.  \eqno (1.3) $$
Next, define the ``maximal correlation coefficient'' [10]
$$ \rho({\cal A}, {\cal B}) :=
\sup |{\rm Corr}(g,h)| \eqno (1.4) $$
where the supremum is taken over all pairs of real-valued, 
square-integrable random variables $g$ and $h$ such
that $g$ is ${\cal A}$-measurable and
$h$ is ${\cal B}$-measurable.
Finally, define
$$ \beta({\cal A}, {\cal B}) :=
\sup {1 \over 2} \sum_{i=1}^I \sum_{j=1}^J
|P(A_i \cap B_j) - P(A_i)P(B_j)|   \eqno (1.5) $$
as well as the ``coefficient of information''
(see e.g.\ [21] or [13])
$$ I({\cal A}, {\cal B}) :=
\sup \sum_{i=1}^I \sum_{j=1}^J
P(A_i \cap B_j) \log \biggl(
{{P(A_i \cap B_j)} \over {P(A_i)P(B_j)}} \biggl)   
\eqno (1.6) $$
where in each of (1.5) and (1.6) the supremum is taken over
all pairs of finite partitions $\{A_1, A_2, \dots, A_I\}$
and $\{B_1, B_2, \dots, B_J\}$ of $\Omega$ such that
$A_i \in {\cal A}$ for each $i$ and 
$B_j \in {\cal B}$ for each $j$.
(Here and below, ``log'' denotes the natural logarithm.)\ \ 
In (1.6) the summand is taken to be 0 if either
$P(A_i)$ or $P(B_j)$ is 0.
It is well known (see e.g.\ [3, v1, Proposition 3.11
and Theorem 5.3(III)]) that
for any two $\sigma$-fields ${\cal A}$ and ${\cal B}$,
$$ \eqalignno{
4\alpha({\cal A}, {\cal B}) &\leq
\rho({\cal A}, {\cal B}), \quad {\rm and} & (1.7) \cr
2\alpha({\cal A}, {\cal B}) &\leq \beta({\cal A}, {\cal B})
\leq \sqrt {I({\cal A}, {\cal B})}. & (1.8) \cr 
}$$

   (D) Now again suppose $X := (X_k, k \in {\bf Z})$ is a 
strictly stationary sequence of ${\bf R}^m$-valued random variables.  (No assumptions on moments.)\ \
For any integer $j$, define the $\sigma$-fields
${\cal F}_{-\infty}^j := \sigma(X_k, k \leq j)$ and
${\cal F}_j^\infty := \sigma(X_k, k \geq j)$.
(Here and below, $\sigma(\dots)$ denotes the $\sigma$-field
$\subset {\cal F}$ generated by $(\dots)$.)\ \ 
For each positive integer $n$, define the following five dependence coefficients:
$$ \eqalignno{ 
\alpha(n) = \alpha(X,n) &:= 
\alpha({\cal F}_{-\infty}^0, {\cal F}_n^\infty); & (1.9) \cr 
\rho(n) = \rho(X,n) &:= 
\rho({\cal F}_{-\infty}^0, {\cal F}_n^\infty); & (1.10) \cr  
\beta(n) = \beta(X,n) &:= 
\beta({\cal F}_{-\infty}^0, {\cal F}_n^\infty); & (1.11) \cr  
I(n) = I(X,n) &:= 
I({\cal F}_{-\infty}^0, {\cal F}_n^\infty); \quad {\rm and} 
& (1.12) \cr
\rho^*(n) = \rho^*(X,n) &:= \sup  
\rho(\sigma(X_k, k \in \Gamma), \sigma(X_k, k \in \Delta)) 
& (1.13) \cr    
}$$
where the supremum in (1.13) is taken over all pairs of
nonempty, disjoint subsets $\Gamma$ and $\Delta$ of ${\bf Z}$ such that 
${\rm dist}(\Gamma,\Delta) := 
\min_{g \in \Gamma, h \in \Delta}|g-h| \geq n$.
(The sets $\Gamma$ and $\Delta$ can be ``interlaced,'' 
i.e.\ with each one containing elements between ones 
in the other set.)\ \ 
Of course by strict stationarity, 
$\alpha(n) = \alpha({\cal F}_{-\infty}^j, 
{\cal F}_{j+n}^\infty)$ for any integer $j$; and the
analogous comment applies to (1.10), (1.11), and (1.12).

The given strictly stationary sequence $X$ is
said to satisfy \hfil\break
``strong mixing'' [23] if 
$\alpha(n) \to 0$ as $n \to \infty$, \hfil\break 
``$\rho$-mixing'' [15] if 
$\rho(n) \to 0$ as $n \to \infty$, \hfil\break
``absolute regularity'' [29] if 
$\beta(n) \to 0$ as $n \to \infty$, \hfil\break
``information regularity'' [21] [29] if 
$I(n) \to 0$ as $n \to \infty$, and \hfil\break
``$\rho^*$-mixing'' [24] [25] if 
$\rho^*(n) \to 0$ as $n \to \infty$. \hfil\break
(The mixing condition in [24] looked somewhat
different from $\rho^*$-mixing, but turned out to be
equivalent to it in the context in that paper;
see [3, v1, Theorem 5.13].)\ \ By (1.7)-(1.8) and
(1.9)-(1.13), the following implications hold: \hfil\break
(i) $\rho^*$-mixing implies $\rho$-mixing, \hfil\break
(ii) $\rho$-mixing implies strong mixing, \hfil\break
(iii) information regularity implies absolute regularity,
and \hfil\break 
(iv) absolute regularity implies strong mixing.
\hfil\break

   With the possible exception of information regularity,
all of these conditions have played a major role in 
limit theory for weakly dependent random variables; see
e.g.\ the books [1], [3], [9], [17], and [22].
Information regularity is sometimes a handy tool in 
the study of stationary Gaussian sequences; see e.g.\ 
[13, Chapter 4] or [3, v3, Chapter 27]. 

    Peligrad [19, Corollary 2.3] proved a central limit
theorem for strictly stationary sequences of real-valued,
square-integrable random variables satisfying the
dependence assumptions
$\rho^*(1) < 1$ and $\alpha(n) \to 0$ as $n \to \infty$.
That result was generalized to strictly stationary
random fields of real-valued random variables 
by Perera [20, Proposition 3] (with the sums being taken 
over a broad class of sets of indices, not just 
``rectangular blocks''). 
It was generalized again in [3, v3, Corollary 29.33]
--- again to strictly stationary random fields of 
real-valued random variables --- 
with another, less restrictive generalization
(to random fields) of the
dependence coefficient $\rho^*(1)$ (but with the sums
taken over just the usual ``rectangular blocks''
of indices).  
Later, for an arbitrary positive integer $m$, 
Tone [26, Theorem 1.1] generalized that latter result to strictly stationary random fields of 
${\bf R}^m$-valued random variables.    
For simplicity, we shall state her result here for
just the special case of random sequences:   
\hfil\break

  {\sc Theorem 1.3} (Tone [26]; Peligrad [19]
for $m = 1$).  
{\sl Suppose $m$ is a positive integer. 
Suppose $X := (X_k, k \in {\bf Z})$ is a strictly 
stationary sequence of ${\bf R}^m$-valued random variables
such that $EX_0 = {\bf 0}_m$ and $E\|X_0\|^2 < \infty$,
and the covariance matrix $\Sigma_{X(0)}$ is 
positive definite (hence nonsingular).
Suppose also that $\rho^*(X,1) < 1$ and that
$\alpha(X,n) \to 0$ as $n \to \infty$.
Then the following two statements hold:

   (1) For each $n \in {\bf N}$, the covariance matrix
$\Sigma_{S(X,n)}$ is positive definite (hence nonsingular).

   (2) One has that (see Notations 1.1(C) and 1.2(A))
$$ \Sigma_{S(X,n)}^{-1/2}S(X,n) \Rightarrow
N({\bf 0}_m, I_m)\ \ {\rm as}\ n \to \infty.  
\eqno (1.14) $$} 

Here in (1.14), the notation $\Rightarrow$ means 
convergence in distribution
on (the Borel $\sigma$-field of) ${\bf R}^m$, and
the notation $N({\bf 0}_m, I_m)$ refers to the 
multivariate normal distribution on ${\bf R}^m$ whose
mean vector is ${\bf 0}_m$ and whose covariance matrix 
is the $m \times m$ identity matrix $I_m$.
The left side of (1.14) is an ${\bf R}^m$-valued random
variable (``random $m \times 1$ column vector'') resulting
from the matrix multiplication indicated there.

   Under different dependence assumptions, again in the
more general context of strictly stationary random fields, 
Bulinskii and Kryzhanovskaya [6, eq.\ (1.13) and Theorem 2]
reformulated a multivariate central limit theorem 
in [7] into the 
form (1.14), with the same use of the ``standard $-1/2$ power'' of the covariance matrix $\Sigma_{S(X,n)}$ as
``normalization,'' and then treated a related central 
limit theorem of the form (1.14) involving the use of the
``standard $-1/2$ power'' of a {\it sample\/}
covariance matrix $\hat \Sigma_{S(X,n)}$ as
``normalization.''
(Those results will not be treated further here.)

   Here is our main result (recall Notations 1.1(E)): 
\hfil\break
   
   {\sc Theorem 1.4.} {\sl Suppose $m$ is a positive 
integer.
Suppose $a$, $b$ and $\tau$ are positive real numbers 
such that $a < b$.
Then there exists a strictly stationary Gaussian
sequence
$X := (X_k, k \in {\bf Z})$ of ${\bf R}^m$-valued,
mean-${\bf 0}_m$
random variables with the following properties:

   (1) $\rho^*(X,1) < 1$.
   
   (2) $\max \{I(X,1), \beta(X,1), \alpha(X,1), \rho(X,1)\}
\leq \tau$. 
     
   (3) $\max \{I(X,n), \beta(X,n), \alpha(X,n), \rho(X,n)\} \to 0$ as $n \to \infty$.

   (4) For every element ($m \times m$ matrix)
$G \in \Lambda_{(m,a,b)}$, there exists an infinite set
$Q \subset {\bf N}$ such that
$$\Sigma_{S(X,n)/\sqrt n} \to G\ \ {\rm as}\ \   
n \to \infty,\ n \in Q. \eqno (1.15) $$}

   Statements (2) and (3) have some redundancy 
(see (1.7)-(1.8)), but that is harmless.
Of course (1.15) means that for every
$(i,j) \in \{1, 2, \dots, m\}^2$, the $(i,j)$-entry
of the matrix $\Sigma_{S(X,n)/\sqrt n}$ converges to
the $(i,j)$-entry of the matrix $G$ as 
$n \to \infty,\ n \in Q$.
Also, the statement that $X$ is a ``Gaussian sequence''
means of course that for any positive integer $L$
and any distinct integers $k(1), k(2), \dots, k(L)$,
the joint distribution of the random vectors
$X_{k(1)}, X_{k(2)}, \dots, Z_{k(L)}$ is a (possibly
degenerate) multivariate normal distribution on
${\bf R}^{Lm}$. 

   Theorem 1.4 will be proved in Section 3, after some
preliminary work is done in Section 2.
In the rest of Section 1 here, a few comments on
this theorem will be given.

   Under the assumptions of Theorem 1.3, 
Tone [26, Claim 3.1] showed that for the covariance 
matrices $\Sigma_{S(X,n)}$, the ratio of the largest to smallest eigenvalues is bounded, and that in fact there 
exists a pair of positive numbers $a < b$ such that
$\Sigma_{S(X,n)/\sqrt n} \in \Lambda_{(m,a,b)}$
for all $n \in {\bf N}$.
Thus in property (4) in Theorem 1.4, the restriction
to matrices in $\Lambda_{(m,a,b)}$ (for some pair of
positive numbers $a < b$) is unavoidable.

   In Theorem 1.4, property (3) cannot be
extended to include $\rho^*(X,n) \to 0$ as 
$n \to \infty$, for that (in conjunction with certain other
properties in Theorem 1.4) would force the covariance
matrices $\Sigma_{S(X,n)/\sqrt n}$ to converge to a
limiting matrix as $n \to \infty$ (a fact
implicitly contained in another, somewhat related result 
of Tone [27, Theorem 3.2]), contradicting property (4).
Also, in Theorem 1.4, the larger the ratio $b/a$ is,
the closer $\rho^*(X,1)$ has to be to 1.
That insight ultimately goes back (in light of basic 
results in [15]) to work of Moore [18]
involving a closely related condition.

   For random sequences and random fields respectively, 
classes of examples constructed in 
[3, v3, Theorem 26.8] and [4, Theorem 1.9]
``separate'' various different but related mixing
assumptions used in [2], [3], [19], [20], [26], 
[27], and other related works.
In particular, the latter class of examples (in [4])
``separates'' the two generalizations 
(to random fields) of the dependence coefficient 
$\rho^*(1)$ (in [20], and in [3] and [26]) implicitly 
alluded to prior to Theorem 1.3.    

   In (1.15), regardless of whether or not the 
eigenvalues of $G$ are simple, 
one can trivially consider a further subsequence
in which the eigenvalues and $m$ orthogonal unit
eigenvectors of the matrices $\Sigma_{S(X,n)/\sqrt n}$
all converge; by a simple calculation, their limits must 
be the eigenvalues and $m$ orthogonal unit eigenvectors
of $G$.
As a consequence, in Theorem 1.3, for the covariance
matrices $\Sigma_{S(X,n)}$, the relative magnitudes of
the eigenvalues, and the respective (orthogonal) 
directions of their eigenvectors, can range essentially
arbitrarily --- within some upper bound (as noted above)
on the ratio of the largest to smallest eigenvalues.
In this respect, Theorem 1.4 helps to ``separate''
Theorem 1.3 from other, more conventional multivariate
central limit theorems (such as the one in 
[27, Theorem 3.2] alluded to above)
in which there is a ``limiting covariance matrix.'' 

   It was noted above that in the special case of
real-valued random variables (i.e.\ $m=1$), Theorem 1.3
boils down to a central limit theorem of
Peligrad [19].
The author [2] (see also [3, v3, Theorem 27.12])
gave a construction (a variant of ones in
[11] and [5]) that showed that in that result of
Peligrad, the growth of the variances need not be 
asymptotically linear, but can instead ``wobble''
between two different linear rates of growth.
That construction was in spirit (though not fully in
letter) a version of Theorem 1.4
for the case $m=1$ (real-valued random variables).

   As was noted above, Theorem 1.3 is actually just a 
special case of a result of Tone [26, Theorem 1.1], which
in its full generality involved random fields 
(of ${\bf R}^m$-valued random variables) indexed
by ${\bf Z}^d$ for an arbitrary positive integer $d$.
By modifying the arguments below, one can prove a
version of Theorem 1.4 for such random fields for
arbitrary ($m$ and) $d$.
However, in the case $d \geq 2$, for such a construction,
the information in Theorem 1.4 that pertains to the
dependence coefficients $\beta(n)$ and $I(n)$
unavoidably becomes false and has to be omitted; see
[3, v3, Theorem 29.9].

   As a simple corollary of Theorem 1.4 itself, one
can derive a version of Theorem 1.4 in which the
sequence $X$ is not Gaussian.
One can simply apply Theorem 1.4 itself with $a$ replaced
by some number $a' \in (0,a)$, then fix $\varepsilon > 0$
such that $a' + \varepsilon < a$, and then replace $X_k$
by $X_k + [V_k^{(1)}, V_k^{(2)}, \dots, V_k^{(m)}]^t$
where $(V_k^{(i)}, k \in {\bf Z}, i \in \{1,\dots,m\})$
is a family of independent, identically distributed
real-valued random variables, this family being
independent of the sequence $X$, with the $V_k^{(i)}$'s
each taking the values $\sqrt \varepsilon$ and
$-\sqrt \varepsilon$ with probability $1/2$ each.
\hfil\break

   {\bf 2. Preliminaries.}  This section will lay some
groundwork for the proof, in Section 3, of Theorem 1.4.  
   
   The random sequence $X$ described
in Theorem 1.4 will be constructed (in Section 3) from
a family of independent ``building block'' random sequences
of a relatively simple structure.
The following lemma will play a role in that 
process of ``assembly.'' \hfil\break

   {\sc Lemma 2.1.} {\sl Suppose 
$(\Omega, {\cal F}, P)$ is a probability space, 
$L$ is a positive integer, and
${\cal A}_\ell$ and ${\cal B}_\ell$,
$\ell \in \{1, 2, \dots, L\}$ are $\sigma$-fields
($\subset {\cal F}$) such that the
$\sigma$-fields ${\cal A}_\ell \vee {\cal B}_\ell$,
$\ell \in \{1, \dots, L\}$ are independent.
Then
$$ \eqalignno{
\rho\Bigl(\bigvee_{\ell = 1}^L {\cal A_\ell},
\bigvee_{\ell = 1}^L {\cal B_\ell}\Bigl)
&= \max_{1 \leq \ell \leq L} \rho({\cal A}_\ell,
{\cal B}_\ell), \quad {\rm and} & (2.1) \cr
I\bigl(\bigvee_{\ell = 1}^L {\cal A_\ell},
\bigvee_{\ell = 1}^L {\cal B_\ell}\Bigl)
&= \sum_{\ell = 1}^L I({\cal A}_\ell,
{\cal B}_\ell). & (2.2) \cr
}$$}

   Proofs of these equalities can be found e.g.\ in
[3, v1, Theorems 6.1 and 6.2(VIII)].
Eq.\ (2.1) is due to Cs\'aki and Fischer
[8, Theorem 6.2].
Eq.\ (2.2) is a classic fact from information theory;
see e.g.\ its role in Pinsker [21].

     The ``building blocks'' for the construction (in
Section 3) of the sequence $X$ for Theorem 1.4 will
be stationary Gaussian sequences of centered
{\it real\/}-valued random variables.
They will be identified (in Section 3) via a careful
choice of their spectral densities.
The rest of Section 2 here will lay some groundwork for
that procedure.
\hfil\break

   {\sc Notations 2.2.} With slight abuse of terminology,
a real Borel function $f$ on
$[-\pi, \pi]$ will be said to be ``symmetric'' if
$f(-\lambda) = f(\lambda)$ for a.e.\ 
$\lambda \in [-\pi, \pi]$.
\smallskip 
   	
   (A) Suppose $f$ is a
real, nonnegative, Borel, symmetric, integrable 
function on $[-\pi, \pi]$.
Suppose $W := (W_k, k \in {\bf Z})$ is a strictly
stationary sequence of real-valued, centered,
square-integrable random variables.
Then $f$ is a ``spectral density function'' for the
sequence $W$ if the following holds:
$$ \forall k \in {\bf Z}, \quad
EW_kW_0 = \int_{-\pi}^{\pi} e^{ik\lambda} f(\lambda) 
{{d\lambda} \over {2\pi}}.  \eqno (2.3) $$
If $W$ has a spectral density function, then it will be
unique modulo sets of Lebesgue measure 0.
The convention on spectral density used here is as in
[3]; it differs by a factor of $2\pi$ from a more
standard convention used in other references.
\smallskip

   (B) For each positive integer $n$, define the real,
nonnegative, symmetric, continuous function 
(the Fej\'er kernel)
$F_n$ on $[-\pi, \pi]$ as follows:
$$ F_n(\lambda) := \cases {
(1/n) \cdot [\sin^2(n \lambda/2)]/
[\sin^2(\lambda/2)] & if 
$\lambda \in [-\pi, \pi] - \{0\}$ \cr
n & if $\lambda = 0$. \cr}  \eqno (2.4) $$

   (C) It is well known that if $W$ and $f$ are as in
(A) above, with $f$ being the spectral density 
function of $W$, then for each positive integer $n$,
$$ E[(W_1 + W_2 + \dots + W_n)/\sqrt n\thinspace]^2 =
\int_{-\pi} ^ \pi F_n(\lambda) f(\lambda) 
{{d\lambda} \over {2\pi}}.
\eqno (2.5) $$
See e.g.\ 
[3, v1, the Note after Lemma 8.18].
\hfil\break

   {\sc Lemma 2.3.} {\sl Suppose $W := (W_k, k \in {\bf Z})$
is a stationary real mean-zero Gaussian random sequence 
that has a spectral density $f$ on $[-\pi, \pi]$
that is bounded a.e.\ between two positive constants.
Then $\rho^*(W,1) < 1$.} 
\hfil\break

   An elementary proof of this lemma can be found in
[3, v1, Theorem 9.8(III)].  (It yields the inequality
$\rho^*(W,1) \leq 1 - a/b$ where $0 < a < b$ and
$a \leq f \leq b$ a.e.  The sharper inequality 
$\rho^*(W,1) \leq (1 - a/b)/(1 + a/b)$
holds as a result of a more sophisticated argument of 
Moore [18] in a closely related context.)

   The analysis that follows will now
involve certain real, Borel, symmetric functions $f$ 
on $[-\pi, \pi]$ that
can take (perhaps even exclusively) {\it negative\/} 
values --- with the intent to use, for some such 
functions $f$ later on, the positive function 
$\lambda \mapsto e^{f(\lambda)}$ 
as the spectral density for a stationary Gaussian
sequence. 
\hfil\break

   {\sc Notations 2.4.} (A) For any 
(not necessarily nonnegative) real, Borel,
square-integrable, symmetric function $f$ on 
$[-\pi, \pi]$, define the quantity
$$ \Psi(f) := \sum_{k=1}^\infty k\psi_{f,k}^2  \eqno (2.6) $$
where for each $k \in {\bf N}$,
$$ \psi_{f,k} := 2 \cdot \int_{-\pi}^\pi 
e^{ik\lambda} f(\lambda) 
{{d\lambda} \over {2\pi}}.  \eqno (2.7) $$
Of course $\sum_{k=1}^\infty \psi_{f,k}^2 < \infty$; 
and with 
$\psi_{f,0} := 
(2 \pi)^{-1} \int_{-\pi}^\pi f(\lambda) d\lambda$, 
one has that \hfil\break 
$\sum_{k=0}^\infty \psi_{f,k} \cos(k \lambda)$
converges in ${\cal L}^2$ to $f$ (and one can say more).
However, the quantity $\Psi(f)$ may be infinite.
\smallskip

   (B) For any two real, Borel, square-integrable, 
symmetric functions $f$ and $g$ on $[-\pi, \pi]$,
one has that $\psi_{f+g,k} = \psi_{f,k} + \psi_{g,k}$
for each $k$ (see (2.7)), and by (2.6) and Minkowski's
inequality,
$[\Psi(f+g)]^{1/2} \leq [\Psi(f)]^{1/2}
+ [\Psi(g)]^{1/2}$ (where if necessary, 
$\infty^{1/2} := \infty$).
\smallskip

   (C) Suppose $a$ and $b$ are real numbers such that
$a < b$.
Suppose $f, f_1, f_2, f_3, \dots$ is a sequence of
real, Borel, symmetric functions on $[-\pi, \pi]$
that are each bounded a.e.\ between $a$ and $b$,
and $f_n \to f$ a.e.\ as $n \to \infty$.
If $\tau$ is a positive number and 
$\Psi(f_n) \leq \tau$ for every $n \in {\bf N}$,
then $\Psi(f) \leq \tau$.  

   (This formulation is unnecessarily restrictive, but will fit our applications later on.
The point is that for each $k$, $\psi_{f(n),k}$ 
(where $f(n)$ means $f_n$) converges to $\psi_{f,k}$
as $n \to \infty$, and hence for each positive integer $L$,
$\sum_{k=1}^L k \psi_{f,k}^2 \leq \tau$, and hence the 
same is true with $L$ replaced by $\infty$.)
\smallskip

   (D) If $f$ is a real, Borel, square-integrable, 
symmetric function on $[-\pi, \pi]$ such that
$\Psi(f) < \infty$, then 
$\int_{-\pi}^\pi e^{f(\lambda)} d\lambda < \infty$.
(This is a special case of a classic result of
Lebedev and Milin [16].
For a detailed exposition of this, see e.g.\ 
[3, v3, Appendix, Theorem A2744(VII)].)
\hfil\break 

   {\sc Lemma 2.5.} {\sl For every $\varepsilon > 0$, there 
exists $\delta = \delta(\varepsilon) > 0$ such that
the following holds:

   Suppose $W:= (W_k, k \in {\bf Z})$ is a stationary
real mean-zero Gaussian random sequence with a 
spectral density function $g$ of the form
$g(\lambda) = e^{f(\lambda)}$, $\lambda \in [-\pi, \pi]$,
where $f$ is a real, Borel, square-integrable, 
symmetric function on $[-\pi, \pi]$ such that
$\Psi(f) \leq \delta$; then 
$I(W,1) \leq \varepsilon$.}  
\hfil\break

   This lemma is implicitly contained in arguments of
Ibragimov, Rozanov, and Solev in [12][14] (see
also [13, Chapter 4]).
A detailed, explicit proof of this lemma can be found in
[3, v3, Theorem 27.11].
\hfil\break 

   {\sc Lemma 2.6.}  {\sl Suppose $\Upsilon_1$, and 
$\Upsilon_2$, and $\theta$ are real numbers such that 
$$ \Upsilon_1 < \theta < \Upsilon_2.   \eqno (2.8)$$ 

   Suppose $\delta > 0$ and $\varepsilon > 0$.
   
   Suppose $N$ is a positive integer.
   
   Suppose $f$ is a real, continuous, symmetric
function on $[-\pi, \pi]$ such that
$$ \eqalignno{
&\Upsilon_1 < f(\lambda) < \Upsilon_2\ 
{\rm for\ all}\ \lambda \in [-\pi, \pi] \quad {\rm and}
& (2.9) \cr
& \Psi(f) < \delta. & (2.10) \cr } $$

   Then there exists a real, continuous, symmetric function
$ h = h_{(f, \Upsilon(1), \Upsilon(2), \theta, \delta,
\varepsilon, N)} $
on $[-\pi, \pi]$ (where the notations $\Upsilon(1)$ and 
$\Upsilon(2)$ mean $\Upsilon_1$ and
$\Upsilon_2$) with the following five properties: 
$$ \eqalignno{
&{\rm For\ every}\ \lambda \in [-\pi, \pi],\ \ 
\Upsilon_1 < h(\lambda) < \Upsilon_2; & (2.11) \cr
&\Psi(h) < \delta; & (2.12) \cr
&|h(0) - \theta| < \varepsilon; & (2.13) \cr
&\int_{-\pi}^{\pi} |h(\lambda) - f(\lambda)| d\lambda < \varepsilon;\ \ {\rm and} & (2.14) \cr
&{\rm for\ every}\ n \in \{1,2,\dots,N\}, \quad  
\Bigl|\int_{-\pi}^{\pi} F_n(\lambda) \cdot e^{h(\lambda)}d\lambda -
\int_{-\pi}^\pi F_n(\lambda) \cdot e^{f(\lambda)}d\lambda
\Bigl|
< \varepsilon. \indent & (2.15)\cr} $$}

   {\sc Proof.}  Refer to (2.8) and (2.9).  
We shall first carry out the proof of
Lemma 2.6 under the following extra assumption: 
$$ \theta > f(0).  \eqno (2.16) $$

   Since $f$ is (by assumption) continuous on the closed
interval $[-\pi, \pi]$, it follows (see (2.8), (2.9), and (2.16)) that there exists a number $c_0$ (henceforth fixed) with the following three properties:
$$ \eqalignno{
& 0 < c_0 < \min\{1, \theta - f(0)\}; & (2.17) \cr
& \Upsilon_1 < f(\lambda) - c_0 < f(\lambda) + c_0
< \Upsilon_2\ \ {\rm for\ all}\ \lambda \in [-\pi, \pi];
\ \ {\rm and}   & (2.18) \cr
&|f(\lambda) - f(0)| < \Upsilon_2 - \theta\ \ 
{\rm for\ all}\ \lambda \in [-c_0, c_0]. & (2.19) \cr
} $$ 
   For each $c \in (0, c_0]$, define the positive numbers
$a_{c,k}, k \in {\bf N}$ as follows:
$$  a_{c,k} := \cases {
(c^2/\pi) \cdot (1/k) & if $k = 1$ or $2$ \cr
(c^2/\pi) \cdot (1/k) \cdot 1/(\log k) & if $k \geq 3$. \cr}
\eqno (2.20) $$
Then for each $c \in (0, c_0]$, one has by (2.17) 
and (2.20) that
$$ \theta - f(0) > a_{c,1} > a_{c,2} > a_{c,3} > \dots
\downarrow 0  \eqno (2.21) $$
and that $\sum_{k=1}^\infty a_{c,k} = \infty$.
Accordingly, for each $c \in (0, c_0]$, let $M(c)$ denote
the greatest positive integer such that (see the first
inequality in (2.21))
$$ \sum_{k=1}^{M(c)} a_{c,k} \leq \theta - f(0).  
\eqno (2.22) $$
For each $c \in (0, c_0]$, define the real, continuous,
symmetric function $g_c$ on $[-\pi, \pi]$ as follows:
For $\lambda \in [-\pi, \pi]$,
$$ g_c(\lambda) := 
\sum_{k=1}^{M(c)} a_{c,k} \cos(k \lambda).  \eqno (2.23) $$ 

   Now suppose $c$ is an arbitrary fixed number such that
$c \in (0, c_0]$.
From (2.23), (2.20), the monotonicity in (2.21), and a standard fact for trigonometric series with nonnegative,
monotonically decreasing coefficients 
(see [3, v3, Appendix, Lemma A2712] --- take the real
parts there --- or [30, p.\ 3, Theorem (2.2)]),
one has that for any $\lambda \in [c, \pi]$,
$$ |g_c(\lambda)| \leq (\pi/\lambda) \cdot a_{c,1}
= (\pi/\lambda) \cdot (c^2 / \pi) \leq c. $$ 
Next suppose for just a moment that $\lambda \in (0,c]$.
Then $0 < \lambda \leq c \leq c_0 < 1$ by (2.17).
Let $I$ denote the positive integer such that
$I < 1/\lambda \leq I+1$.  
Then for all $k \in \{1, 2, \dots, I\}$, one has that
$k\lambda < 1$ and hence $\cos(k \lambda) > 0$.
If $M(c) \leq I$, then it follows from (2.23) and (2.20)
that $g_c(\lambda) > 0$.
If instead $M(c) > I$, then one has
$\sum_{k=1}^I a_{c,k} \cos (k \lambda) > 0$ and
(since $1 \leq \lambda \cdot (I+1)$) again by (2.20),
(2.17), and the monotonicity in (2.21),
(again see [3, v3, Lemma A2712] or [30, p.\ 3])
$$ \Bigl| \sum_{k=I+1}^{M(c)} 
a_{c,k} \cos (k \lambda) \Bigl|
\leq (\pi/\lambda) \cdot a_{c, I+1}
\leq (\pi/\lambda) \cdot (c^2/\pi) \cdot
(1/(I+1)) \leq c^2 < c,$$
and hence $g_c(\lambda) \geq -c$ by (2.23).
Putting all these pieces together (see also (2.22)
and (2.23) again), one now has that
$$ \eqalignno{
&|g_c(\lambda)| \leq c\ \ {\rm for\ all}\ 
\lambda \in [c, \pi];\ \ {\rm and} & (2.24) \cr
&-c \leq g_c(\lambda) \leq \sum_{k=1}^{M(c)}a_{c,k}
\leq \theta - f(0)\ \ {\rm for\ all}\ \lambda \in [0,c].
& (2.25) \cr
} $$ 
(Eq.\ (2.25) was shown above for $\lambda \in (0,c]$; it
extends to $\lambda = 0$ by continuity of the function
$g_c$.)\ \ 
By (2.18), (2.24), and (2.18) again (keeping in mind our
ongoing assumption $c \in (0, c_0]$), one has that for
all $\lambda \in [c, \pi]$,
$$ \Upsilon_1 < f(\lambda) -c \leq f(\lambda) + g_c(\lambda)
\leq f(\lambda) + c < \Upsilon_2.$$
By (2.18), (2.25), and (2.19), for all $\lambda \in [0,c]$,
$$ \Upsilon_1 < f(\lambda) -c \leq f(\lambda) + g_c(\lambda)
 \leq f(\lambda) + \theta - f(0) 
 < \Upsilon_2 - \theta + \theta = \Upsilon_2. $$
Hence by symmetry, one now has that
$$ \Upsilon_1 < f(\lambda) + g_c(\lambda) < \Upsilon_2\ \ 
{\rm for\ all}\ \lambda \in [-\pi, \pi].   
\eqno (2.26) $$

   Equations (2.24), (2.25), and (2.26) were shown for
any arbitrary $c \in (0, c_0]$.
Our plan now is to let the function $h$ be defined by
$$h := f + g_c   \eqno (2.27) $$
for some sufficiently small $c \in (0, c_0]$.
To start off, note that under (2.27) for any given 
$c \in (0, c_0]$, (2.11) holds by (2.26).

   Next, by (2.20), for each $c \in (0, c_0]$, 
$$ \sum_{k=1}^\infty k \cdot a_{c,k}^2 =
(c^4/\pi^2) \cdot
\Bigl[ 1 + (1/2) + 
\sum_{k=3}^\infty 1/[k(\log k)^2] \Bigl] 
< \infty; $$   
and in fact the middle term converges to 0 as $c \to 0+$.
Hence by (2.23) and (2.6)-(2.7), $\Psi(g_c) \to 0$
as $c \to 0+$.  
Hence by (2.10) and Notations 2.4(B), 
$[\Psi(f + g_c)]^{1/2} < \delta^{1/2}$ for all
$c \in (0, c_0]$ sufficiently small.
Thus under (2.27), eq.\ (2.12) holds for all
$c \in (0, c_0]$ sufficiently small.

   Next, for each $c \in (0, c_0]$, by the definition
of the positive integer $M(c)$ (see the entire sentence
containing (2.22)), followed by (2.20), one has that
$$ 0 \leq [\theta - f(0)] - \sum_{k=1}^{M(c)} a_{c,k}
< \sum_{k=1}^{M(c)+1} a_{c,k} - \sum_{k=1}^{M(c)} a_{c,k}
= a_{c,M(c)+1} \leq c^2/\pi. $$
That is, by (2.23),
$ 0 \leq [\theta-f(0)] - g_c(0) < c^2/\pi$, that is,
$ 0 \leq \theta - [f(0) + g_c(0)] < c^2/\pi.$
Hence under (2.27), eq.\ (2.13) holds for all
$c \in (0, c_0]$ sufficiently small.

   Next, by (2.24) and symmetry, for every 
$\lambda \in [-\pi, \pi] - \{0\}$,
$g_c(\lambda) \to 0$ as $c \to 0+$.
Hence by (2.9), (2.26), and dominated convergence,
(2.14) holds (under (2.27)) for all $c \in (0, c_0]$
sufficiently small.
Also, since each Fej\'er kernel (see (2.4)) is bounded,
and by (2.9) and (2.26) the functions $\exp(f(\lambda))$
and $\exp(f(\lambda) + g_c(\lambda))$ 
(for $c \in (0, c_0]$) 
are uniformly bounded (between 
$\exp \Upsilon_1$ and $\exp \Upsilon_2$),
one has by dominated convergence that (under (2.27))
eq.\ (2.15) holds for all $c \in (0, c_0]$
sufficiently small.
Thus under (2.27), eqs.\ (2.11)--(2.15)
hold for all $c \in (0, c_0]$ sufficiently small.
Thus Lemma 2.6 holds under the extra assumption (2.16).

   It will be useful to note that, again under the
extra assumption (2.16), one can expand the statement of
Lemma 2.6 to include the following variant of (2.15):
$$ {\rm For\ every}\ n \in \{1,2,\dots,N\}, 
\quad  
\Bigl|\int_{-\pi}^{\pi} F_n(\lambda) \cdot e^{-h(\lambda)}d\lambda -
\int_{-\pi}^\pi F_n(\lambda) \cdot e^{-f(\lambda)}d\lambda
\Bigl|
< \varepsilon. \indent \eqno (2.28) $$
To accomplish this, one shows that under (2.27), 
eq.\ (2.28) holds for all 
$c \in (0, c_0]$ sufficiently small.
The argument is essentially the same as the corresponding
one for (2.15) in the preceding paragraph. 

   Now let us briefly take care of the cases where (2.16)
does not hold.
Refer to (2.8) and (2.9) again.
If $\theta = f(0)$, then let $h := f$ and we are done.
Finally, if $\theta < f(0)$, then by replacing 
$\Upsilon_1$, $\Upsilon_2$, $\theta$, and $f$ by
$-\Upsilon_2$, $-\Upsilon_1$, $-\theta$, and $-f$
(note that $\Psi(-f) = \Psi(f)$ by (2.6)-(2.7)), 
one trivially converts to the case where (2.16) holds.
(The resulting function, say $\tilde h$, is then
multiplied by $-1$ to produce the final function $h$.
In order for (2.15) to result at the end of this
``trivial conversion argument,'' it was vital to
derive the ``extra'' fact (2.28) at the end of the
argument under (2.16) above.)\ \  
That completes the proof of Lemma 2.6.  
\bigskip

   {\bf 3. Proof of Theorem 1.4.} The proof will be
written out here for the case $m \geq 2$.
(The argument for the case $m = 1$ is similar but less
complicated.)\ \  
The proof will be divided into several ``steps.'' 
(One of those ``steps'' will be a ``lemma.'')
\hfil\break

   {\sc Step 3.1.} Refer to the statement of Theorem 1.4.
Decreasing $\tau$ and/or $a$ and/or increasing $b$ if necessary, we assume without loss of generality that
$$ \eqalignno{
&0 < a < 1 < b \quad {\rm and} & (3.1) \cr
&0 < \tau < 1.  & (3.2) \cr
}$$

   Let us identify the set of all (real) $m \times m$ 
matrices with ${\bf R}^{m \uparrow 2}$ 
(with each entry in the matrix identified with a 
coordinate in ${\bf R}^{m \uparrow 2}$).
The set ${\bf R}^{m \uparrow 2}$ is separable.
Hence every nonempty subset of  ${\bf R}^{m \uparrow 2}$
is separable (an elementary fact --- see e.g.
[3, v3, Appendix, Lemma A3101]).
Accordingly, let $\tilde \Lambda$ be a countable dense
subset of $\Lambda_{(m,a,b)}$.
Let $G_1, G_2, G_3, \dots$ be a sequence of elements of
$\tilde \Lambda$ such that (for convenience) each
element of $\tilde \Lambda$ is listed infinitely many
times in that sequence.

   In order to prove Theorem 1.4, it suffices to construct a 
strictly stationary, mean-${\bf 0}_m$ Gaussian sequence
$X := (X_k, k \in {\bf Z})$ of ${\bf R}^m$-valued
random variables such that properties (1), (2), and (3)
in Theorem 1.4 hold as well as the following property:
($4'$) There exists a strictly increasing sequence 
$(N_1, N_2, N_3, \dots)$ of positive integers, and
a positive number $\Theta$, such that
(recall Notations 1.1(F))  
for all $n \in {\bf N}$ sufficiently large, 
$$\Sigma_{S(X,N(n))/\sqrt{N(n)}} - G_n \in 
{\bf B}_{\rm sym}^{(m)}[2^{-n}\Theta].   
\eqno (3.3) $$
(Here and throughout the rest of this note, when the
notation $N_n$ appears in a subscript, it will be written
$N(n)$ for typographical convenience.)\ \  
It will then follow trivially that  
each member $G \in \tilde \Lambda$ would be the
limit of a subsequence of the matrices
$\Sigma_{S(X,N(n))/\sqrt{N(n)}}$ (for the integers $n$ such that $G_n = G$); and property (4) in Theorem 1.4 would then
follow as an easy consequence.

   We shall return to the matrices $G_n$ in Step 3.5 below.
\hfil\break

   {\sc Step 3.2.} Refer again to (3.1).
In what follows, for convenience, our attention will be ``expanded'' from $\Lambda_{(m,a,b)}$ to to
$\Lambda_{(m,a/2,2b)}$.

   Define the positive number
$$ \gamma := a/(20m^2).  \eqno (3.4) $$

   Define the (``lattice'') set
$$ {\bf L} := \{\dots, -3\gamma, -2\gamma, -\gamma, 0,
\gamma, 2\gamma, 3\gamma, \dots \} \eqno (3.5) $$
(that is, the set of all real numbers of the form $k\gamma$,
$k \in {\bf Z}$).
Let $\Lambda_{\bf L}$ denote the set of all $m \times m$
matrices
$H := (h_{ij}, 1 \leq i,j \leq m) \in \Lambda_{(m,a/2,2b)}$
such that $h_{ij} \in {\bf L}$ for every
$(i,j) \in \{1, \dots, m\}^2$.
By Notations 1.1(D)(E) (see the third sentence after
(1.2)), the set $\Lambda_{(m,a/2,2b)}$ is bounded
(as represented as a subset of ${\bf R}^{m \uparrow 2}$).
It follows that $\Lambda_{\bf L}$ is a finite set.
Of course the set $\Lambda_{\bf L}$ is nonempty.
(For example, $cI_m \in \Lambda_{\bf L}$ where $c$ is
an element of ${\bf L}$ such that $a/2 < c < a$ ---
such a $c$ exists by (3.4).)\ \ 
Define the positive integer
$$ L := {\rm card}\thinspace \Lambda_{\bf L}.  \eqno (3.6) $$
Let the elements of $\Lambda_{\bf L}$ be denoted as
$Q_1^{(1)}, Q_2^{(1)}, \dots, Q_L^{(1)}$,
with the representation
$$ Q_{\ell}^{(1)} := 
(q^{(1)}_{\ell i j}, 1 \leq i,j \leq m)  \eqno (3.7) $$
for $\ell \in \{1,2,\dots, L\}$.
These matrices $Q^{(1)}_\ell$ are of course 
symmetric and positive definite (since they belong to 
$\Lambda_{(m,a/2,2b)}$). 
\hfil\break

   {\sc Step 3.3.}  Two other classes of matrices will be
needed.  (These matrices will be symmetric but not
positive definite.)

   For each $u \in \{1, \dots, m\}$, let
$Q_u^{(2)} := (q_{uij}^{(2)}, 1 \leq i,j \leq m)$ denote
the (symmetric) $m \times m$ matrix defined by
$$ q_{uij}^{(2)} := \cases{
1 & if $(i,j) = (u,u)$ \cr
0 & for all other $(i,j)$. \cr}
\eqno (3.8) $$

   Now recall the assumption $m \geq 2$ made in the 
first sentence of Section 3.
Let ${\bf T}$ denote the set of all ordered pairs 
$(u,v) \in \{1,\dots,m\}^2$ such that $u < v$.
For each ordered pair $(u,v) \in {\bf T}$, let
$Q_{uv}^{(3)} := (q_{uvij}^{(3)}, 1 \leq i,j \leq m)$ 
denote the (symmetric) $m \times m$ matrix defined by
$$ q_{uvij}^{(3)} := \cases{
1 & if $(i,j) \in \{(u,u), (u,v), (v,u), (v,v)\}$ \cr
0 & for all other $(i,j)$. \cr}
\eqno (3.9) $$

   Now to set the stage for the next lemma (and for some
other calculations below), note that trivially by (3.1) 
and (3.4),
$\gamma/(2bL) < \gamma < 10m\gamma < 1$.
\hfil\break 
   
   {\sc Lemma 3.4.} {\sl For every matrix 
$G \in \Lambda_{(m,a,b)}$, there exists an array
$$ c = c(G) := 
\Bigl\{ c_{\ell}^{(1)}, \ell \in \{1,2,\dots, L\}; \thinspace
c_u^{(2)}, u \in \{1,2, \dots, m\}; \thinspace
c_{uv}^{(3)}, (u,v) \in {\bf T} \Bigl\} \eqno (3.10) $$
of positive numbers such that the following statements
hold:
$$ \eqalignno{
& \forall \ell \in \{1,\dots, L\},\ \ 
\gamma/(2bL) \leq c_\ell^{(1)} \leq 1; & (3.11) \cr
& \forall u \in \{1, \dots, m\},\ \  
2m\gamma \leq c_u^{(2)} \leq 10m\gamma; & (3.12) \cr
&\forall (u,v) \in {\bf T},\ \ 
2\gamma \leq c_{uv}^{(3)} \leq 5\gamma;  
& (3.13) \cr
}$$  
and
$$ G = \sum_{\ell = 1}^L c_{\ell}^{(1)} Q_{\ell}^{(1)} 
 + \sum_{u = 1}^m c_u^{(2)} Q_u^{(2)}
 + \sum_{(u,v) \in {\bf T}} c_{uv}^{(3)} Q_{uv}^{(3)}.
 \eqno (3.14) $$}
 
   {\sc Proof.} Represent the matrix $G$ by
$$ G := (g_{ij}, 1 \leq i,j \leq m).  \eqno (3.15) $$
Of course by the hypothesis and Notations 1.1(E),
$G$ is symmetric.   
For each $(i,j) \in \{1,\dots,m\}^2$, let $\kappa_{ij}$
denote the integer such that (see (3.4))
$$ \kappa_{ij}\gamma \leq g_{ij} < (\kappa_{ij}+1)\gamma.
\eqno (3.16) $$
Then $\kappa_{ij} = \kappa_{ji}$.
Define the (symmetric) $m \times m$ matrix
$H := (h_{ij}, 1 \leq i,j \leq m)$ as follows:
$$ \eqalignno{ 
&\forall i \in \{1,\dots,m\},\ \ 
h_{ii} := (\kappa_{ii} -8m)\gamma; \quad {\rm and} 
& (3.17) \cr
&\forall (i,j) \in {\bf T},\ \ 
h_{ij} = h_{ji} := (\kappa_{ij} -3)\gamma. & (3.18) \cr
}$$    
Now for each $(i,j) \in \{1,\dots,m\}^2$,
$$ |g_{ij} - h_{ij}| \leq
|g_{ij} - \kappa_{ij}\gamma| 
+ |\kappa_{ij}\gamma - h_{ij}|.
\eqno (3.19) $$ 
In the right side of (3.19), the first term is bounded 
above by $\gamma$ (by (3.16)), and the second term is
either $8m\gamma$ (if $i=j$) or  $3 \gamma$ (if $i \neq j$),
by (3.17)-(3.18).
Hence $G-H \in {\bf B}_{\rm sym}^{(m)}[8m\gamma]$.
Recall from (3.4) and (3.1) that $8m^2\gamma < a/2 < b$.
Since (by hypothesis) $G \in \Lambda_{(m,a,b)}$, 
it now follows from Notations 1.1(G) that 
$H \in \Lambda_{(m,a/2,2b)}$.
Hence by (3.17)-(3.18) and the sentence after (3.5),
$H \in \Lambda_{\bf L}$.
Accordingly (see the sentence after (3.6)) let
$\ell'$ denote the element of $\{1,\dots,L\}$ such that
$$ Q_{\ell'}^{(1)} = H.  \eqno (3.20) $$

   Define the array $c = c(G)$ in (3.10) 
(in a slightly unconventional order) as follows:
First,
$$ c_{\ell'}^{(1)} := 1 \quad {\rm and} \quad
\forall \ell \in \{1,\dots, L\} - \{\ell'\},\ \ 
c_\ell^{(1)} := \gamma /(2bL).  \eqno (3.21) $$
Next, for convenience, referring to (3.7), 
define the $m \times m$ symmetric
matrix $S := (s_{ij}, 1 \leq i,j \leq m)$ 
as follows:
$$ \forall (i,j) \in \{1, \dots, m\}^2,\ \ 
s_{ij} := \sum_{\ell \in \{1,\dots, L\} - \{\ell'\}}
c_\ell^{(1)}q_{\ell ij}^{(1)}.  \eqno (3.22) $$
(By (3.21), 
$S = [\gamma/(2bL)]
\sum_{\ell \in \{1,\dots, L\} - \{\ell'\}}
Q_\ell^{(1)}$;
however, the form (3.22) will be a little more natural
for the calculations that follow.)\ \  
Next, use $S$ to continue the definition of the array in (3.10) as follows:
$$ \forall (u,v) \in {\bf T},\ \ 
c_{uv}^{(3)} := [g_{uv} - h_{uv}] - s_{uv}.  \eqno (3.23) $$
Finally, use (3.23) itself to complete the definition of 
the array in (3.10) as follows: 
$$ \forall u \in \{1,\dots,m\},\ \ 
c_u^{(2)} := [g_{uu} - h_{uu}] - s_{uu} 
- \sum_{(i,j) \in {\bf T}: u \in \{i,j\}} c_{ij}^{(3)}.
\eqno (3.24) $$

   Now recall from the entire last paragraph of Step 3.2
that $Q_\ell \in \Lambda_{(m,a/2,2b)}$ for every
$\ell \in \{1, \dots, L\}$.
It follows from (3.7) and Notations 1.1(D)(E) (see the 
third sentence after (1.2)) that for each 
$\ell \in \{1, \dots, L\}$
and each $(i,j) \in \{1, \dots, m\}^2$,
$|q_{\ell ij}^{(1)}| \leq 2b$.
Hence by (3.21) and (3.22),
$$ \forall (i,j) \in \{1, \dots, m\}^2,\ \ 
|s_{ij}| \leq \gamma, \eqno (3.25) $$
that is, $S \in {\bf B}_{\rm sym}^{(m)}[\gamma]$.

   Now we shall verify eqs.\ (3.11)-(3.14) (though not
quite in that order).

   First, (3.11) holds by (3.21) and the sentence 
after (3.9).

   Next, for each $(u,v) \in {\bf T}$, by (3.23) and (3.18),
$$ c_{uv}^{(3)} =
[g_{uv} - \kappa_{uv}\gamma] 
+ [\kappa_{uv} \gamma - h_{uv}] - s_{uv} 
= [g_{uv} - \kappa_{uv}\gamma] + 3\gamma -s_{uv}. 
\eqno (3.26) $$
By (3.16) and (3.25), the far right side of (3.26) is
bounded below by $0 + 3\gamma - \gamma$ and bounded
above by $\gamma + 3\gamma + \gamma$.
Hence (3.13) holds.

   Next let us verify (3.12). 
For any given $u \in \{1, \dots ,m\}$, the following holds:  
The set $\{(i,j) \in {\bf T}: u \in \{i,j\}\}$ has
exactly $m-1$ elements ($(1,u), \dots, (u-1,u)$ and
$(u,u+1), \dots, (u,m)$), and hence by (3.13) (just proved 
above),
$$ 2(m-1)\gamma 
\leq \sum_{(i,j) \in {\bf T}: u \in \{i,j\}} c_{ij}^{(3)}
\leq 5(m-1)\gamma.  \eqno (3.27) $$
Now by (3.24),
$$ c_u^{(2)} = [g_{uu} - \kappa_{uu}\gamma]
+ [\kappa_{uu}\gamma - h_{uu}] - s_{uu}
- \sum_{(i,j) \in {\bf T}: u \in \{i,j\}} c_{ij}^{(3)}. 
\eqno (3.28) $$
By (3.16), (3.17), (3.25), and (3.27), the right side 
of (3.28) is bounded below by 
$0 + 8m\gamma - \gamma -5m\gamma$
and bounded above by 
$\gamma + 8m\gamma + \gamma - 0$.
Hence (3.12) holds.

   Finally, (3.14) needs to be verified.
First, for $(i,j) \in {\bf T}$, 
by (3.21), (3.22), (3.8), (3.9), (3.20) (with (3.7)),
and (3.23),  
$$ \eqalignno{
\sum_{\ell = 1}^L &c_{\ell}^{(1)} q_{\ell ij}^{(1)} 
 + \sum_{u = 1}^m c_u^{(2)} q_{uij}^{(2)}
 + \sum_{(u,v) \in {\bf T}} c_{uv}^{(3)} q_{uvij}^{(3)}
 =
 1 \cdot q_{\ell' ij}^{(1)} + s_{ij} + 0 + 
 c_{ij}^{(3)} \cdot 1 \cr 
&= h_{ij} + s_{ij} + c_{ij}^{(3)} = g_{ij}. & (3.29) \cr 
 }$$
Next, recall that the matrices $G$, $Q_\ell^{(1)}$, 
$Q_u^{(2)}$, and $Q_{uv}^{(3)}$ (and $S$) are symmetric. 
As a trivial consequence, for $(i,j) \in {\bf T}$,
the far left and far right sides of (3.29) remain equal 
if the indices $i$ and $j$ are switched.
Finally, for each $i \in \{1,\dots, m\}$,
by (3.21) (again with (3.20)), (3.22), (3.8), (3.9), 
and (3.24),  
$$
\sum_{\ell = 1}^L c_{\ell}^{(1)} q_{\ell ii}^{(1)} 
 + \sum_{u = 1}^m c_u^{(2)} q_{uii}^{(2)}
 + \sum_{(u,v) \in {\bf T}} c_{uv}^{(3)} q_{uvii}^{(3)}  
= 1 \cdot h_{ii} + s_{ii} + c_i^{(2)} \cdot 1
+ \sum_{(u,v) \in {\bf T}: i \in \{u,v\}}
c_{uv}^{(3)} \cdot 1
= g_{ii}.
$$
From all of these observations, (3.14) holds.
That completes the proof of Lemma 3.4.
\hfil\break

   {\sc Step 3.5.}  This step will involve, after some
preliminary work, repeated
applications of Lemma 2.6.  The notation
$h_{(f, \Upsilon(1), \Upsilon(2), \theta, \delta,
\varepsilon, N)}$ in Lemma 2.6 (see the sentence after
(2.10)) will be used repeatedly, and for typographical
convenience it will be written below as
$h(f, \Upsilon_1, \Upsilon_2, \theta, \delta,
\varepsilon, N)$.

   For the use of that notation, define
(see (3.1), (3.2), (3.4), and (3.6)) the real numbers
$$ \Upsilon_1 := \log \Bigl(\gamma/(3bL)\Bigl); \quad
\Upsilon_2 := \log 2; \quad
{\rm and} \quad 
\delta := \delta\Bigl(\tau^2/[2m(L + 1 + m)]\Bigl) 
\eqno (3.30) $$
where in the last equality we are using the notation
in Lemma 2.5.
By (3.1), (3.4), and (3.30), 
$\Upsilon_1 < 0 < \Upsilon_2$.
Referring to (2.6)-(2.7), 
we shall say that a given real, continuous, symmetric 
function $f$ on $[-\pi, \pi]$ satisfies
``Condition {\bf C}'' if (2.9) and (2.10) hold for the
given values in (3.30).

   Next, refer to Notations 2.2(B), involving the
Fej\'er kernels.
Of course by Fej\'er's Theorem, if $f$ is a (say) real, continuous, symmetric function on $[-\pi, \pi]$, then
$(2\pi)^{-1} \int_{-\pi}^{\pi} F_n(\lambda) \cdot
\allowbreak f(\lambda) d\lambda$
converges to $f(0)$ as $n \to \infty$.
For a given real, continuous, symmetric function
$f$ on $[-\pi, \pi]$ and a given $\varepsilon > 0$,
let ${\cal N}(f,\varepsilon)$ be a positive integer such that
$$ \forall n \geq {\cal N}(f,\varepsilon), \quad
\Bigl| f(0) - \int_{-\pi}^{\pi} F_n(\lambda) 
f(\lambda) {{d\lambda} \over {2\pi}} 
\Bigl| \leq \varepsilon. 
\eqno (3.31) $$
   
   Next, refer to the sequence $G_1, G_2, G_3, \dots$ of
matrices in $\Lambda_{(m,a,b)}$ (in fact in 
$\tilde \Lambda$) from the second paragraph of
Step 3.1. 
Applying Lemma 3.4 and using the notations there,
define for each positive integer $n$ the 
array 
$$ {\bf c}_n = c(G_n) := 
\Bigl\{ c_{\ell,n}^{(1)}, \ell \in \{1,2,\dots, L\}; \thinspace
c_{u,n}^{(2)}, u \in \{1,2, \dots, m\}; \thinspace
c_{u,v,n}^{(3)}, (u,v) \in {\bf T} \Bigl\} \eqno (3.32) $$
of positive numbers (satisfying (3.11)-(3.14)
with $G = G_n$).
By (3.11), (3.12), and (3.13), together with (3.30) and 
the sentence after (3.9), one has that for each positive integer $n$ and
each number $c$ in the array ${\bf c}_n$, 
$\Upsilon_1 < \log c \leq 0 < \Upsilon_2$. 
  
   Now we shall define a sequence of positive integers
$(N_0, N_1, N_2, \dots)$; and 
we shall define, for each positive 
integer $n$, a collection 
$${\cal C}_n := \Bigl\{ f_{\ell,n}^{(1)}, \ell \in \{1,2,\dots, L\};
f_{u,n}^{(2)}, u \in \{1,2,\dots,m\};
f_{u,v,n}^{(3)}, (u,v) \in {\bf T} \Bigl\} 
\eqno (3.33) $$
of real, continuous, symmetric functions on $[-\pi, \pi]$
that each satisfy Condition {\bf C}
(see the sentence after (3.30)).
Notice that for a given positive integer $n$, there will 
be only finitely many functions in this array (3.33) ---
in fact $L + m + m(m-1)/2$ of them.
The definition will be recursive in $n$, with
$N_{n-1}$ and ${\cal C}_n$ being defined together
for $n = 1, 2, 3, \dots\thinspace $.
It proceeds as follows:

   To start off, define the positive integer 
$N_0 := 1$,
and let each of the functions in the collection 
${\cal C}_1$ in (3.33) be
the trivial constant function with range $\{0\}$.
Of course a constant function $f$ on $[-\pi, \pi]$
satisfies $\Psi(f) = 0$.
Since $\Upsilon_1 < 0 < \Upsilon_2$ (as was noted above), 
it now follows that the (constant) functions 
in (3.33) (for $n=1$) satisfy Condition {\bf C}. 

   Now suppose $n \geq 1$ is an integer, and the positive
integer $N_{n-1}$ and the real, continuous, symmetric
functions in ${\cal C}_n$ in (3.33) have already been
defined, and that those functions all satisfy
Condition {\bf C}.
Define the positive integer
$$N_n := N_{n-1} + \max {\cal N}(e^f,2^{-n})   
\eqno (3.34)$$
where this maximum is taken over all functions 
$f$ in the collection ${\cal C}_n$ in (3.33) for 
the given $n$.
(Of course for each such $f$, the notation $e^f$ simply
refers to the real, continuous, symmetric function
$\lambda \mapsto e^{f(\lambda)}$ on $[-\pi,\pi]$.)
Now referring to (3.30), (3.32), and the sentence 
after (3.32), and 
applying Lemma 2.6, define the functions in 
the collection ${\cal C}_{n+1}$ as follows:
First, for each $\ell \in \{1,\dots,L\}$, 
define the function $f_{\ell,n+1}^{(1)}$ by
$$ f_{\ell,n+1}^{(1)} :=
h(f_{\ell,n}^{(1)}, \Upsilon_1, \Upsilon_2, 
\log c_{\ell,n+1}^{(1)}, \delta, 2^{-n}, N_n). 
\eqno (3.35) $$    
Next, for each $u \in \{1, \dots, m\}$, define the
function $f_{u,n+1}^{(2)}$ by
$$ f_{u,n+1}^{(2)} :=
h(f_{u,n}^{(2)}, \Upsilon_1, \Upsilon_2, 
\log c_{u,n+1}^{(2)}, \delta, 2^{-n}, N_n). 
\eqno (3.36) $$
Finally, for each $(u,v) \in {\bf T}$, define the
function $f_{u,v,n+1}^{(3)}$ by
$$ f_{u,v,n+1}^{(3)} :=
h(f_{u,v,n}^{(3)}, \Upsilon_1, \Upsilon_2, 
\log c_{u,v,n+1}^{(3)}, \delta, 2^{-n}, N_n). 
\eqno (3.37) $$
That completes the definition of the collection
${\cal C}_{n+1}$.  Note that from 
(2.11)-(2.12) in Lemma 2.6, each
of the functions in this collection ${\cal C}_{n+1}$
satisfies Condition {\bf C}.
   
   That completes the recursive definition of the
positive integers $N_0, N_1, N_2, \dots$ and the 
collections ${\cal C}_n, n \in {\bf N}$.
From (3.34) and the definition of $N_0$, one
has that
$$ 1 = N_0 < N_1 < N_2 < \dots\thinspace . \eqno (3.38) $$

   {\sc Step 3.6.} The next task is to establish a
collection 
$${\cal C} := \Bigl\{ f_\ell^{(1)}, \ell \in 
\{1,2,\dots, L\};
f_u^{(2)}, u \in \{1,2,\dots,m\};
f_{u,v}^{(3)}, (u,v) \in {\bf T} \Bigl\} 
\eqno (3.39) $$
of ``limit functions'' on $[-\pi, \pi]$ from the 
collections ${\cal C}_n$.

   First suppose $\ell \in \{1, \dots, L\}$.
For each positive integer $n$, from (3.35) and 
eq.\ (2.14) in Lemma 2.6, one has that
$\int_{-\pi}^\pi 
|f_{\ell,n+1}^{(1)}(\lambda) - f_{\ell,n}^{(1)}(\lambda)|   
d\lambda < 2^{-n}$.
Hence
$\int_{-\pi}^\pi \sum_{n=1}^\infty 
|f_{\ell,n+1}^{(1)}(\lambda) - f_{\ell,n}^{(1)}(\lambda)|
d\lambda < \infty$.
Hence $\sum_{n=1}^\infty 
|f_{\ell,n+1}^{(1)}(\lambda) - f_{\ell,n}^{(1)}(\lambda)|
< \infty$ for a.e.\ $\lambda \in [-\pi, \pi]$.
Define the function $f_{\ell}^{(1)}$ 
a.e.\ on $[-\pi, \pi]$ as follows:
$$ f_{\ell}^{(1)}(\lambda) := \lim_{n \to \infty}   
f_{\ell,n}^{(1)}(\lambda). \eqno (3.40) $$
The right side of (3.40) will be defined in ${\bf R}$
for a.e.\ $\lambda \in [-\pi, \pi]$.  
On the null-set of values $\lambda$ for which that
limit does not exist in ${\bf R}$, the quantity
$f_{\ell}^{(1)}(\lambda)$ is left undefined here.

   Next, for each $u \in \{1,\dots, m\}$, going through
the same procedure, but using (3.36) instead of (3.35),
define the function $f_u^{(2)}$ a.e.\ on $[-\pi, \pi]$ by
$$ f_u^{(2)}(\lambda) := \lim_{n \to \infty}   
f_{u,n}^{(2)}(\lambda). \eqno (3.41) $$

   Finally, for each $(u,v) \in {\bf T}$,  
again going through the same procedure,
this time using (3.37),
define the function $f_{u,v}^{(3)}$ a.e.\ on 
$[-\pi, \pi]$ by
$$ f_{u,v}^{(3)}(\lambda) := \lim_{n \to \infty}   
f_{u,v,n}^{(3)}(\lambda). \eqno (3.42) $$

   That completes the definition of the collection
${\cal C}$ in (3.39).
Since each of the functions in each of the collections 
${\cal C}_n$
is real and symmetric and satisfies Condition {\bf C}, it follows from (3.40)-(3.42) that each of the functions 
in the collection ${\cal C}$ is
a.e.\ real and symmetric, with its range being bounded a.e.\ 
within the closed interval $[\Upsilon_1, \Upsilon_2]$.
\hfil\break

   {\sc Step 3.7.}  Next some calculations involving
Fej\'er kernels will be given.
Later on, they will play a key role in obtaining bounds
on the covariance matrices for partial sums of sequences
of random vectors (${\bf R}^m$-valued random variables).

   For each positive integer $n$, define the array
$$ {\bf c}^*_n := 
\Bigl\{ c_{\ell,n}^{*(1)}, \ell \in \{1,2,\dots, L\}; \thinspace
c_{u,n}^{*(2)}, u \in \{1,2, \dots, m\}; \thinspace
c_{u,v,n}^{*(3)}, (u,v) \in {\bf T} \Bigl\} \eqno (3.43) $$
of positive numbers as follows:
First, for each $\ell \in \{1, \dots, L\}$,
referring to (3.38), (3.40), and (2.4), define the 
positive number
$$ c_{\ell,n}^{*(1)} := 
\int_{-\pi}^\pi F_{N(n)}(\lambda)\cdot 
\exp (f_\ell^{(1)}(\lambda))
{{d\lambda} \over {2\pi}}.    \eqno (3.44) $$
Next, for each $u \in \{1, \dots, m\}$, referring to (3.41),
define the positive number  
$$ c_{u,n}^{*(2)} := 
\int_{-\pi}^\pi F_{N(n)}(\lambda)\cdot 
\exp (f_u^{(2)}(\lambda))
{{d\lambda} \over {2\pi}}.    \eqno (3.45) $$
Finally, for each $(u,v) \in {\bf T}$, referring to (3.42),
define the positive number  
$$ c_{u,v,n}^{*(3)} := 
\int_{-\pi}^\pi F_{N(n)}(\lambda)\cdot 
\exp (f_{u,v}^{(3)}(\lambda))
{{d\lambda} \over {2\pi}}.   \eqno (3.46) $$
That completes the definition of the array ${\bf c}_n^*$ 
in (3.43). 

   Our next task is to compare the arrays ${\bf c}_n$
and ${\bf c}_n^*$ in (3.32) and (3.43).

   To start that process, suppose $n \geq 2$, and  
suppose $\ell \in \{1, \dots, L\}$.
By (3.34), 
$N_n > {\cal N}(\exp f_{\ell,n}^{(1)}, 2^{-n})$;
and hence by (3.31),
$$ \Bigl|\exp (f_{\ell,n}^{(1)}(0)) 
- \int_{-\pi}^\pi F_{N(n)}(\lambda) \cdot
\exp(f_{\ell,n}^{(1)}(\lambda))
{{d\lambda} \over {2\pi}} 
\Bigl| \leq 2^{-n}.  \eqno (3.47)$$
Also, for each integer $p \geq n$, one has that
$f_{\ell,p+1}^{(1)} :=
h(f_{\ell,p}^{(1)}, \Upsilon_1, \Upsilon_2, 
\log c_{\ell,p+1}^{(1)}, \delta, 2^{-p}, N_p)$
by (3.35), and since $N_n \leq N_p$ by (3.38)
one therefore has from eq.\ (2.15) in Lemma 2.6 that
$$ \Bigl| \int_{-\pi}^\pi F_{N(n)} \cdot
\exp (f_{\ell,p+1}^{(1)}(\lambda))
{{d\lambda} \over {2\pi}}
- \int_{-\pi}^\pi F_{N(n)} \cdot
\exp (f_{\ell,p}^{(1)}(\lambda))
{{d\lambda} \over {2\pi}} \Bigl| \leq 2^{-p}. 
\eqno (3.48) $$
By (3.47) and (3.48), using a telescoping sum, 
one has that  
$$ \eqalignno{ \forall p \geq n+1, \quad
\Bigl|\exp (f_{\ell,n}^{(1)}(0)) 
&- \int_{-\pi}^\pi F_{N(n)}(\lambda) \cdot
\exp(f_{\ell,p}^{(1)}(\lambda))
{{d\lambda} \over {2\pi}} 
\Bigl| \cr 
&\leq 2^{-n} + [2^{-n} + 2^{-(n+1)} + \dots + 2^{-(p-1)}]
\leq 3 \cdot 2^{-n}. \indent \indent  
& (3.49) \cr}$$
Now recall that for each $p \geq 1$, the function 
$f_{\ell, p}^{(1)}$ satisfies Condition {\bf C} and
is therefore bounded between $\Upsilon_1$ and 
$\Upsilon_2$, and hence the function
$\exp f_{\ell, p}^{(1)}$ is bounded between
$\exp \Upsilon_1$ and $\exp \Upsilon_2$. 
Since any given Fej\'er Kernel is bounded,
one now has by (3.44), (3.40), (3.49), and dominated 
convergence (taking the limit as $p \to \infty$) that
for our given fixed $n$ and $\ell$,
$$
|\exp (f_{\ell,n}^{(1)}(0)) - c_{\ell,n}^{*(1)}| = 
\Bigl|\exp (f_{\ell,n}^{(1)}(0)) 
- \int_{-\pi}^\pi F_{N(n)}(\lambda) \cdot
\exp(f_{\ell}^{(1)}(\lambda))
{{d\lambda} \over {2\pi}} 
\Bigl|  
\leq 3 \cdot 2^{-n}. \indent \eqno (3.50) 
$$
Recall our supposition here that $n \geq 2$.
From (3.35) (with $n+1$ replaced by $n$)
and eq.\ (2.13) in Lemma 2.6, one has that
$$|f_{\ell,n}^{(1)}(0) - \log c_{\ell,n}^{(1)}| 
< 2^{-(n-1)}. 
\eqno (3.51) $$  
Since $f_{\ell,n}^{(1)}$ satisfies condition {\bf C},
one trivially has (see (3.30) 
and the sentence after it)
that $f_{\ell,n}^{(1)}(0) < \log 2$.
From the sentence after (3.32), one also
has that $\log c_{\ell,n}^{(1)} \leq 0 < \log 2$.
Since $de^x/dx = e^x \leq 2$ for $x \leq \log 2$,
it now follows from (3.51) and trivial calculus that
$|\exp(f_{\ell,n}^{(1)}(0)) - c_{\ell,n}^{(1)}| 
\leq 4 \cdot 2^{-n}.$
Hence by (3.50), 
$|c_{\ell,n}^{*(1)} - c_{\ell,n}^{(1)}| 
\leq 7 \cdot 2^{-n}$.

   Let us display for convenient reference what we
have just verified:
$$ \forall n \geq 2,\ \forall \ell \in \{1, \dots, L\},
\quad  
|c_{\ell,n}^{*(1)} - c_{\ell,n}^{(1)}| 
\leq 7 \cdot 2^{-n}.    \eqno (3.52) $$
With arguments exactly analogous to that of (3.52),
using (3.41)-(3.42) and (3.45)-(3.46) in place of
(3.40) and (3.44),
one has that
$$ \forall n \geq 2,\ \forall u \in \{1, \dots, m\},
\quad  
|c_{u,n}^{*(2)} - c_{u,n}^{(2)}| 
\leq 7 \cdot 2^{-n}    \eqno (3.53) $$  
and that
$$ \forall n \geq 2,\ \forall (u,v) \in {\bf T}, \quad  
|c_{u,v,n}^{*(3)} - c_{u,v,n}^{(3)}| 
\leq 7 \cdot 2^{-n}.    \eqno (3.54) $$  
We shall return to (3.52)-(3.54) later on.
\hfil\break

   {\sc Step 3.8.}  Our task in this step is to
construct the random sequence $X$ for Theorem 1.4.  
That will be done with a family of  ``building blocks'' 
that are independent of each other, each one being a 
stationary real mean-zero Gaussian random sequence 
with a particular spectral density function.
We shall use the well known fact that any real,
nonnegative, Borel, symmetric, integrable function on
$[-\pi, \pi]$ is the spectral density function of
some stationary real mean-zero Gaussian random sequence.

   Refer to (3.40), (3.41), and (3.42).
For each $\ell \in \{1, \dots, L\}$ and each
$p \in \{1, \dots, m\}$, let
$X^{(1,\ell,p)} := (X^{(1,\ell,p)}_k, k \in {\bf Z})$
be a stationary real mean-zero Gaussian random 
sequence with spectral density
function $\exp f^{(1)}_\ell$ on $[-\pi, \pi]$.
For each $u \in \{1, \dots, m\}$, let
$X^{(2,u)} := (X^{(2,u)}_k, k \in {\bf Z})$
be a stationary real mean-zero Gaussian random 
sequence with spectral density
function $\exp f^{(2)}_u$ on $[-\pi, \pi]$. 
For each $(u,v) \in {\bf T}$, let
$X^{(3,u,v)} := (X^{(3,u,v)}_k, k \in {\bf Z})$
be a stationary real mean-zero Gaussian random 
sequence with spectral density
function $\exp f^{(3)}_{u,v}$ on $[-\pi, \pi]$.
Let these random sequences be constructed in such a
way that they are all independent of each other.

   Refer to Notations 1.1(C).  For each 
$\ell \in \{1, \dots, L\}$, let
$(Q^{(1)}_\ell)^{1/2}$ denote the symmetric positive
definite $m \times m$ ``square root'' matrix of
$Q^{(1)}_\ell$.
(Recall the sentence after (3.7).)

   Define the sequence $X := (X_k, k \in {\bf Z})$
of ${\bf R}^m$-valued random variables as follows:
For each $k \in {\bf Z}$,

$$ \eqalignno{ 
X_k :=
\sum_{\ell = 1}^L
(Q^{(1)}_\ell)^{1/2}
&\Bigl[X_k^{(1,\ell, 1)}, X_k^{(1,\ell, 2)}, \dots,
X_k^{(1,\ell, m)}\Bigl]^t \cr 
&+ \sum_{u=1}^m X_k^{(2,u)} {\bf e}_u
+ \sum_{(u,v) \in {\bf T}} X_k^{(3,u,v)} 
({\bf e}_u + {\bf e}_v).
& (3.55) \cr
}$$
Here and below, for $p \in \{1, \dots, m\}$, 
${\bf e}_p := [0, \dots, 0, 1, 0, \dots, 0]^t$
where the $1$ is the $p^{\rm th}$ coordinate.
The first sum in the right side of (3.55) involves
matrix multiplication; the other two involve
simple scalar multiplication.    
By elementary arguments, $X$ is a strictly
stationary, Gaussian sequence of ${\bf R}^m$-valued, 
mean-${\bf 0}_m$ random variables. 
Our task now is to verify properties (1)-(4)
stipulated in Theorem 1.4.
The ``mixing properties'' (1)-(3) will be verified in
Step 3.9, and property (4) will be verified in 
Step 3.10.
\hfil\break

   {\sc Step 3.9.}  In this step, the mixing properties
(1), (2), and (3) stipulated in Theorem 1.4 will be
verified (though not in that order).

   For each positive integer $n$, by (3.55) and 
Lemma 2.1 (and the independence of the ``building block''
sequences in the second paragraph of Step 3.8),
$$ \rho(X,n) \leq \max \Bigl\{
\max_{1 \leq \ell \leq L, 1 \leq p \leq m}
\rho(X^{(1,\ell,p)},n), 
\max_{1 \leq u\leq m} \rho(X^{(2,u)},n),
\max_{(u,v)\in {\bf T}} \rho(X^{(3,u,v)},n) \Bigl\}  
\eqno (3.56) $$
and
$$ I(X,n) \leq
\sum_{\ell = 1}^L \sum_{p=1}^m I(X^{(1,\ell,p)},n)
+ \sum_{u=1}^m I(X^{(2,u)},n)
+ \sum_{(u,v) \in {\bf T}} I(X^{(3,u,v)},n).
\eqno (3.57) $$

   Next some calculations connected with information
regularity are needed for the ``building block''
sequences.

    Referring to the sentence containing (3.33), 
one has that for 
each $\ell \in \{1, \dots, L\}$ and each positive 
integer $n$,
$\Upsilon_1 < f^{(1)}_{\ell,n}(\lambda) < \Upsilon_2$ 
for all $\lambda \in [-\pi, \pi]$, and also
$\Psi(f^{(1)}_{\ell, n}) < \delta$. 
Hence for each $\ell \in \{1, \dots, L\}$,  
by (3.40) and Notations 2.4(C), 
$\Psi(f^{(1)}_\ell) \leq \delta$. 
Hence for a given $\ell \in \{1, \dots, L\}$ and a given
$p \in \{1, \dots, m\}$, by the second paragraph
in Step 3.8, one has from (3.30) and Lemma 2.5 that
(i) $I(X^{(1,\ell,p)},1) \leq \tau^2/[2m(L + 1 + m)]$,
and hence by [3, v3, Lemma 27.9(I)(II)],
one also has that 
(ii) $\rho(X^{(1,\ell,p)},1) 
\leq [2I(X^{(1,\ell,p)},1)]^{1/2}
\leq \tau$,
and that (iii) $I(X^{(1,\ell,p)},n) \to 0$ and
$\rho(X^{(1,\ell,p)},n) \to 0$ as $n \to \infty$.
By exactly analogous arguments,
using (3.41) and (3.42) in place of (3.40), 
one obtains (i), (ii), and (iii) with
$X^{(2,u)}$ (for $u \in \{1, \dots, m\}$) and with 
$X^{(3,u,v)}$ (for $(u,v) \in {\bf T}$) in place of
$X^{(1,\ell,p)}$.

   Hence by (3.55) and (3.57) (and (3.2)), 
$$I(X,1) \leq 
[mL + m + m(m-1)/2] \cdot \tau^2/[2m(L + 1 + m)] 
\leq \tau^2 \leq \tau; $$
and hence by (1.8), $\alpha(X,1) \leq \beta(X,1) \leq \tau$; 
and also by (3.55) and (3.56), $\rho(X,1) \leq \tau$.
Also, by (3.55), (3.56), and (3.57), 
$I(X,n) \to 0$ and $\rho(X,n) \to 0$ as $n \to \infty$;
and hence also by (1.8), $\alpha(X,n) \to 0$ and
$\beta(X,n) \to 0$ as $n \to \infty$.
Thus properties (2) and (3) in Theorem 1.4 hold.

   Next, recall from above that for a given
$\ell \in \{1, \dots, L\}$ and a 
given positive integer $n$, one has that
$\Upsilon_1 < f^{(1)}_{\ell,n}(\lambda) < \Upsilon_2$ 
for all $\lambda \in [-\pi, \pi]$. 
Hence by (3.40), for a given $\ell \in \{1, \dots, L\}$,  
$\Upsilon_1 \leq f^{(1)}_\ell (\lambda) \leq \Upsilon_2$
for a.e.\ $\lambda \in [-\pi, \pi]$.
Hence by the second paragraph of Step 3.8, 
for a given $\ell \in \{1, \dots, L\}$ and a given
$p \in \{1, \dots, m\}$, the 
stationary Gaussian sequence $X^{(1, \ell, p)}$ has
a spectral density function that is bounded a.e.\ 
between the two positive constants $\exp \Upsilon_1$
and $\exp \Upsilon_2$,
and hence by Lemma 2.3 it satisfies 
$\rho^*(X^{(1, \ell, p)},1) < 1$.
By exactly analogous arguments, using
(3.41) and (3.42) in place of (3.40), one has that
$\rho^*(X^{(2,u)},1) < 1$ for $u \in \{1, \dots, m\}$
and that    
$\rho^*(X^{(3,u,v)},1) < 1$ for $(u,v) \in {\bf T}$.
Now by (3.55) and Lemma 2.1, eq.\ (3.56) holds
with each $\rho$ replaced by $\rho^*$.
It now follows that $\rho^*(X,1) < 1$.
Thus property (1) in Theorem 1.4 holds.
\hfil\break

   {\sc Step 3.10.}  In this final step, we shall
verify property (4) in Theorem 1.4, by showing that
for the sequence $(N_1, N_2, N_3, \dots)$ of positive
integers defined in Step 3.5 (see (3.38)), there
exists a positive number $\Theta$ such that (3.3) holds
for all $n \geq 2$.

   Refer again to the second paragraph of Step 3.8, where 
the sequences $X^{(1,\ell,p)}$, $X^{(2,u)}$, and 
$X^{(3,u,v)}$ are defined.
One of course has that for each $\ell \in \{1,\dots, L\}$, 
each $p \in \{1, \dots, m\}$, and each $n \in {\bf N}$,
$E[n^{-1/2} S(X^{(1,\ell,p)}, n)] = 0$; and the
analogous comment applies with $X^{(1, \ell, p)}$
replaced by $X^{(2,u)}$ or $X^{(3,u,v)}$.
(That should be kept in mind in the calculations that
follow.)\ \ 
By (3.44) and Notations 2.2(B)(C) (and the second paragraph 
of Step 3.8), for each $\ell \in \{1, \dots, L\}$,
each $p \in \{1, \dots, m\}$ and each $n \in {\bf N}$,
$$ E[N_n^{-1/2} S(X^{(1,\ell,p)}, N_n)]^2
= \int_{-\pi}^\pi 
F_{N(n)}(\lambda) \cdot \exp(f^{(1)}_\ell(\lambda))
{{d\lambda} \over {2\pi}}
= c^{*(1)}_{\ell,n}.  
\eqno (3.58) $$
By similar arguments using (3.45) and (3.46) in place of
(3.44), one has that for each $u \in \{1, \dots, m\}$
and each $n \in {\bf N}$,
$$ E[N_n^{-1/2} S(X^{(2,u)}, N_n)]^2
= \int_{-\pi}^\pi 
F_{N(n)}(\lambda) \cdot \exp(f^{(2)}_u(\lambda))
{{d\lambda} \over {2\pi}}
= c^{*(2)}_{u,n},  
\eqno (3.59) $$
and that for each $(u,v) \in {\bf T}$ and each 
$n \in {\bf N}$,
$$ E[N_n^{-1/2} S(X^{(3,u,v)}, N_n)]^2
= \int_{-\pi}^\pi 
F_{N(n)}(\lambda) \cdot \exp(f^{(3)}_{u,v}(\lambda))
{{d\lambda} \over {2\pi}}
= c^{*(3)}_{u,v,n}.  
\eqno (3.60) $$
 
   In what will now follow, we shall repeatedly use 
the fact that if $V$ is an ${\bf R}^m$-valued random 
variable such that $EV = {\bf 0}_m$ and 
$E\|V\|^2 < \infty$, then the $m \times m$ covariance 
matrix $\Sigma_V$ can be written simply as 
$\Sigma_V = EVV^t$.
   
   For each $\ell \in \{1, \dots, L\}$ and each
$n \in {\bf N}$, define the ${\bf R}^m$-valued random
variable
$$ Y_n^{(\ell)} :=
n^{-1/2}[S(X^{(1,\ell,1)},n), S(X^{(1,\ell,2)},n),
   \dots, S(X^{(1,\ell,m)},n)]^t. $$
By (3.58) and the independence of the sequences
$X^{(1,\ell,p)}$, $p \in \{1, \dots, m\}$ (again see the
second paragraph of Step 3.8), one has that for each
$\ell \in \{1, \dots, L\}$ and each $n \in {\bf N}$,
the ${\bf R}^m$-valued random variable $Y_{N(n)}^{(\ell)}$
has mean vector ${\bf 0}_m$ and covariance matrix
$EY_{N(n)}^{(\ell)}(Y_{N(n)}^{(\ell)})^t
= c^{*(1)}_{\ell,n}I_m$.
Hence for each $\ell \in \{1, \dots, L\}$ and each
$n \in {\bf N}$ (recall that the matrix 
$(Q_\ell^{(1)})^{1/2}$ is symmetric),  
the ${\bf R}^m$-valued random vector
$$ N_n^{-1/2} \sum_{k=1}^{N(n)} 
(Q_\ell^{(1)})^{1/2}
[X_k^{(1,\ell,1)}, X_k^{(1,\ell,2)}, \dots, 
X_k^{(1,\ell,m)}]^t
= (Q_\ell^{(1)})^{1/2}Y_{N(n)}^{(\ell)}$$
has mean vector ${\bf 0}_m$ and covariance matrix
$$ E(Q_\ell^{(1)})^{1/2} Y_{N(n)}^{(\ell)}(Y_{N(n)}^{(\ell)})^t  ((Q_\ell^{(1)})^{1/2})^t  
= (Q_\ell^{(1)})^{1/2} c^{*(1)}_{\ell,n}I_m (Q_\ell^{(1)})^{1/2}
= c^{*(1)}_{\ell,n}Q_\ell^{(1)}.
\eqno (3.61) $$

   By (3.59) and the entire sentence containing (3.8),
for each $u \in \{1, \dots, m\}$, the 
${\bf R}^m$-valued random variable
$N_n^{-1/2} \sum_{k=1}^{N(n)} X_k^{(2,u)} {\bf e}_u
= N_n^{-1/2}[S(X^{(2,u)}, N_n)]{\bf e}_u$
trivially has mean vector ${\bf 0}_m$ and covariance 
matrix
$$ c^{*(2)}_{u,n} {\bf e}_u {\bf e}_u^t 
= c^{*(2)}_{u,n}Q_u^{(2)}.  \eqno (3.62) $$ 
Similarly, by (3.60) and the entire sentence 
containing (3.9), for each $(u,v) \in {\bf T}$,
the ${\bf R}^m$-valued random variable
$N_n^{-1/2} \sum_{k=1}^{N(n)} X_k^{(3,u,v)} 
({\bf e}_u + {\bf e}_v)
= N_n^{-1/2}[S(X^{(3,u,v)}, N_n)]({\bf e}_u + {\bf e}_v)$
has mean vector ${\bf 0}_m$ and covariance matrix
$$ c^{*(3)}_{u,v,n} ({\bf e}_u + {\bf e}_v) 
({\bf e}_u + {\bf e}_v)^t
= c^{*(3)}_{u,v,n}Q_{u,v}^{(3)}.  \eqno (3.63) $$ 

   Now we use the elementary equality
$\Sigma_{Y+Z+ \dots + V} = 
\Sigma_Y + \Sigma_Z + \dots + \Sigma_V$ for an
arbitrary finite collection $Y, Z, \dots, V$ of
independent ${\bf R}^m$-valued random variables whose coordinates have finite second moments.      
By (3.55) and the independence of the sequences in the
second paragraph of Step 3.8, followed by the entire
sentences containing (3.61), (3.62), and (3.63), 
one has that for each
$n \in {\bf N}$, the ${\bf R}^m$-valued random variable
$N_n^{-1/2}S(X,N_n)$ has mean vector ${\bf 0}_m$ and
covariance matrix 
$$ G_n^* := 
\sum_{\ell = 1}^L c^{*(1)}_{\ell,n} Q^{(1)}_\ell 
+ \sum_{u = 1}^m c^{*(2)}_{u,n} Q^{(2)}_u 
+ \sum_{(u,v) \in {\bf T}} c^{*(3)}_{u,v,n} Q^{(3)}_{u,v}.
\eqno (3.64)   
$$
      
   Now from (3.32) and (3.14), for each $n \in {\bf N}$,
$$ G_n = 
\sum_{\ell = 1}^L c^{(1)}_{\ell,n} Q^{(1)}_\ell 
+ \sum_{u = 1}^m c^{(2)}_{u,n} Q^{(2)}_u 
+ \sum_{(u,v) \in {\bf T}} c^{(3)}_{u,v,n} Q^{(3)}_{u,v}.
\eqno (3.65)   
$$

   Recall from the final paragraph of Step 3.2 that
$Q^{(1)}_\ell \in \Lambda_{(m,a/2,2b)}$ for each
$\ell \in \{1, \dots, L\}$.
By (3.7) and Notations 1.1(D)(E) (see the third
sentence after (1.2)), one has that 
$|q_{\ell i j}| \leq 2b$ for each
$\ell \in \{1, \dots, L\}$ and each 
$(i,j) \in \{1, \dots, m\}^2$.
Taking that together with the entire sentences
containing (3.8) and (3.9), and then using
(3.52), (3.53), and (3.54), one obtains from
(3.64) and (3.65) that for each $n \geq 2$,
$$ G_n^* - G_n \in 
{\bf B}^{(m)}_{\rm sym}
\Bigl[\Bigl((2b \cdot L) + (1 \cdot m) + (1 \cdot m(m-1)/2)
\Bigl) \cdot
7 \cdot 2^{-n}\Bigl]. $$
Referring again to the entire sentence 
containing (3.64), one has that for the positive
number $\Theta := 7 \cdot [2bL + m + m(m-1)/2]$, eq.\ (3.3)
holds for all $n \geq 2$.
That completes the proof of property (4) in Theorem 1.4.
The proof of Theorem 1.4 is complete.
\hfil\break 

\centerline {REFERENCES}
\bigskip

\refs [1] H.C.P.\ Berbee, 
{\it Random Walks with Stationary Increments and Renewal Theory\/},
(Mathematical Centre, Amsterdam, 1979).

\refs [2] R.C.\ Bradley, 
On the growth of variances in a central limit theorem 
for strongly mixing sequences,
{\it Bernoulli\/} 
{\bf 5} (1999), 67-80.

\refs [3] R.C.\ Bradley, 
{\it Introduction to Strong Mixing Conditions\/}, 
Volumes 1, 2, and 3,
(Kendrick Press, Heber City, Utah, 2007).

\refs [4] R.C.\ Bradley, 
On the dependence coefficients associated with three 
mixing conditions for random fields,
In:\ {\it Dependence in Probability, Analysis and Number
Theory\/}, 
(I.\ Berkes, R.C.\ Bradley, H.\ Dehling,
M.\ Peligrad, and R.\ Tichy, eds.), 
pp.\ 89-121,
(Kendrick Press, Heber City, Utah, 2010).

\refs [5] W.\ Bryc and A.\ Dembo, 
On large deviations of empirical measures for stationary Gaussian processes,
{\it Stochastic Process.\ Appl.\/}
{\bf 58} (1995), 23-34.

\refs [6] A.\ Bulinskii and N.\ Kryzhanovskaya,
Convergence rate in CLT for vector-valued random fields
with self-normalization,
{\it Probab.\ Math.\ Statist.}
{\bf 26} (2006), 261-281.

\refs [7] A.\ Bulinskii and A.\ Shashkin,
{\it Limit Theorems for Associated Random Fields and
Related Systems\/}, 
Advanced Series on Statistical Science and Applied 
Probability, 10, 
(World Scientific, Hackensack, NJ, 2007).

\refs [8] P.\ Cs\'aki and J.\ Fischer, 
On the general notion of maximal correlation,
{\it Magyar Tud.\ Akad.\ Mat.\ Kutato Int.\ Kozl.\/}
{\bf 8} (1963), 27-51.

\refs [9] P.\ Doukhan, 
{\it Mixing:\ Properties and Examples}, 
(Springer-Verlag, New York, 1995).

\refs [10] H.O.\ Hirschfeld, 
A connection between correlation and contingency,
{\it Proc.\ Camb.\ Phil.\ Soc.\/}
{\bf 31} (1935), 520-524.

\refs [11] I.A.\ Ibragimov,
On the spectrum of stationary Gaussian sequences 
satisfying the strong mixing condition II.\ Sufficient conditions.\ Mixing rate,
{\it Theor.\ Probab.\ Appl.\/}
{\bf 15} (1970), 23-36.

\refs [12] I.A.\ Ibragimov and Yu.A.\ Rozanov,
On the connection between two characteristics of 
dependence of Gaussian random vectors,
{\it Theor.\ Probab.\ Appl.\/}
{\bf 15} (1970), 295-299.

\refs [13] I.A.\ Ibragimov and Yu.A.\ Rozanov,
{\it Gaussian Random Processes\/},
(Springer-Verlag, New York, 1978).

\refs [14] I.A.\ Ibragimov and V.N.\ Solev,
A condition for regularity of a Gaussian stationary
sequence,
{\it Soviet Math.\ Dokl.\/}
{\bf 10} (1969), 371-375.

\refs [15] A.N.\ Kolmogorov and Yu.A. Rozanov,
On strong mixing conditions for stationary Gaussian
processes,
{\it Theor.\ Probab.\ Appl.\/}
{\bf 5} (1960), 204-208.

\refs [16] N.A.\ Lebedev and I.M.\ Milin,
An inequality,
{\it Vestnik Leningrad Univ.\/} 
{\bf 20} (1965), 157-158.

\refs [17] Z.\ Lin and C.\ Lu,
{\it Limit Theory for Mixing Dependent Random 
Variables\/},
(Science Press, Beijing, and 
Kluwer Academic Publishers, Boston, 1996).

\refs [18] C.C.\ Moore,
The degree of randomness in a stationary time series,
{\it Ann.\ Math.\ Statist.\/}
{\bf 34} (1963), 1253-1258.

\refs [19] M.\ Peligrad,
On the asymptotic normality of weak dependent random
variables,
{\it J.\ Theor.\ Probab.\/}
{\bf 9} (1996), 703-715.

\refs [20] G.\ Perera,
Geometry of ${\bf Z}^d$ and a central limit theorem
for weakly dependent random fields,
{\it J.\ Theor.\ Probab.\/}
{\bf 10} (1997), 581-603.

\refs [21] M.S.\ Pinsker,
{\it Information and Information Stability of Random
Variables and Processes\/},
(Nauka, Moscow, 1960).  (In Russian)

\refs [22] E.\ Rio,
{\it Theorie asymptotiques des processus al\'eatoires
faiblements d\'ependantes\/}, \hfil\break
Math\'ematiques \& Applications 31,
(Springer, Berlin, 2000).

\refs [23] M.\ Rosenblatt, 
A central limit theorem and a strong mixing condition,
{\it Proc.\ Natl.\ Acad.\ Sci.\ USA\/}
{\bf 42} (1956), 43-47.

\refs [24] M.\ Rosenblatt, 
Central limit theorems for stationary processes,
{\it Proceedings of the Sixth Berkeley Symposium on
Probability and Statistics\/}, Volume 2,
pp.\ 551-561,
(University of California Press, Los Angeles, 1972).

\refs [25] C.\ Stein,
A bound for the error in the normal approximation to the
distribution of a sum of dependent random variables,
{\it Proceedings of the Sixth Berkeley Symposium on
Probability and Statistics\/}, Volume 2,
pp.\ 583-602,
(University of California Press, Los Angeles, 1972).

\refs [26] C.\ Tone,
A central limit theorem for multivariate strongly
mixing random fields,
{\it Probab.\ Math.\ Statist.\/}
{\bf 30} (2010), 215-222.

\refs [27] C.\ Tone,
Central limit theorems for Hilbert-space valued random
fields satisfying a strong mixing condition,
{\it Lat.\ Am.\ J.\ Prob.\ Math.\ Stat.\/}
{\bf 8} (2011), 77-94.

\refs [28] S.A.\ Utev and M.\ Peligrad,
Maximal inequalities and an invariance principle for a
class of weakly dependent random variables,
{\it J.\ Theor.\ Probab.\/}
{\bf 16} (2003), 101-115.

\refs [29] V.A.\ Volkonskii and Yu.A.\ Rozanov,
Some limit theorems for random functions I,
{\it Theor.\ Probab.\ Appl.\/}
{\bf 4} (1959), 178-197.

\refs [30] A.\ Zygmund,
{\it Trigonometric Series\/}, Volumes 1 and 2,
(Cambridge University Press, Cambridge, 1959).

\bye